\documentstyle[amstex,12pt]{article}

\font\tenfrak=eufm10
\font\sevenfrak=eufm7
\font\fivefrak=eufm5
\newfam\frakfam \def\frak{\fam\frakfam\tenfrak} \textfont\frakfam=\tenfrak
  \scriptfont\frakfam=\sevenfrak  \scriptscriptfont\frakfam=\fivefrak

\begin{document}
    \pagestyle{plain}
    \setlength{\baselineskip}{1.2\baselineskip}
    \setlength{\parindent}{\parindent}

\title{{\bf On the variety of Lagrangian subalgebras}}
\author{Sam Evens\thanks{Research partially supported by NSF 
grants DMS-9623322 and DMS-9970102;}\\
Dept. of Math., University of Notre Dame, Notre Dame, IN 46656\\
and Dept. of Math., University of Arizona, Tucson, AZ 85721\\
\hspace{.04in}
 and Jiang-Hua Lu\thanks{Research partially supported by
NSF grant DMS 9803624.}\\
Dept. of Math., University of Arizona, Tucson, AZ 85721\\
evens{\@}math.arizona.edu, jhlu{\@}math.arizona.edu}
\maketitle

\newtheorem{thm}{Theorem}[section]
\newtheorem{lem}[thm]{Lemma}
\newtheorem{prop}[thm]{Proposition}
\newtheorem{cor}[thm]{Corollary}
\newtheorem{rem}[thm]{Remark}
\newtheorem{exam}[thm]{Example}
\newtheorem{nota}[thm]{Notation}
\newtheorem{dfn}[thm]{Definition}
\newtheorem{ques}[thm]{Question}
\newtheorem{eq}{thm}

\newcommand{\rw}{\rightarrow}
\newcommand{\lrw}{\longrightarrow}
\newcommand{\rhu}{\rightharpoonup}
\newcommand{\lhu}{\leftharpoonup}
\newcommand{\Map}{\longmapsto}
\newcommand{\qed}{\begin{flushright} {\bf Q.E.D.}\ \ \ \ \
                  \end{flushright} }
\newcommand{\beqa}{\begin{eqnarray*}}
\newcommand{\eeqa}{\end{eqnarray*}}

\newcommand{\ZR}{Z_{\Bbb R}}
\newcommand{\fls}{\fl_{d,\sigma}}
\newcommand{\flstriv}{\fl_\sigma}
\newcommand{\flsC}{\fl_{d,\sigma, {\Bbb C}}}
\newcommand{\flsnull}{{\fl}_{d,{\sigma}_0}}
\newcommand{\flsnullC}{{\fl}_{d,{\sigma}_0, {\Bbb C}}}
\newcommand{\frds}{{\fr}_{d,\sigma}}
\newcommand{\Rds}{{R}_{d,\sigma}}
\newcommand{\os}{{\cal O}_{\Bbb C}}
\newcommand{\cO}{{\cal O}}

\newcommand{\la}{\mbox{$\langle$}}
\newcommand{\ra}{\mbox{$\rangle$}}
\newcommand{\ot}{\mbox{$\otimes$}}
\newcommand{\xa}{\mbox{$x_{(1)}$}}
\newcommand{\xb}{\mbox{$x_{(2)}$}}
\newcommand{\xc}{\mbox{$x_{(3)}$}}
\newcommand{\ya}{\mbox{$y_{(1)}$}}
\newcommand{\yb}{\mbox{$y_{(2)}$}}
\newcommand{\yc}{\mbox{$y_{(3)}$}}
\newcommand{\yd}{\mbox{$y_{(4)}$}}
\renewcommand{\aa}{\mbox{$a_{(1)}$}}
\newcommand{\ab}{\mbox{$a_{(2)}$}}
\newcommand{\ac}{\mbox{$a_{(3)}$}}
\newcommand{\ad}{\mbox{$a_{(4)}$}}
\newcommand{\ba}{\mbox{$b_{(1)}$}}
\newcommand{\bt}{\mbox{$b_{(2)}$}}
\newcommand{\bc}{\mbox{$b_{(3)}$}}
\newcommand{\ca}{\mbox{$c_{(1)}$}}
\newcommand{\cb}{\mbox{$c_{(2)}$}}
\newcommand{\cc}{\mbox{$c_{(3)}$}}
\newcommand{\uo}{\mbox{$\underbar{o}$}}
\newcommand{\eo}{\mbox{$e$}}

\newcommand{\ts}{\mbox{$\sigma$}}
\newcommand{\las}{\mbox{${}_{\sigma}\!A$}}
\newcommand{\lasone}{\mbox{${}_{\sigma'}\!A$}}
\newcommand{\ras}{\mbox{$A_{\sigma}$}}
\newcommand{\rds}{\mbox{$\cdot_{\sigma}$}}
\newcommand{\lds}{\mbox{${}_{\sigma}\!\cdot$}}

\newcommand{\bb}{\mbox{$\bar{\beta}$}}
\newcommand{\bg}{\mbox{$\bar{\gamma}$}}

\newcommand{\id}{\mbox{${\rm id}$}}
\newcommand{\Fun}{\mbox{${\rm Fun}$}}
\newcommand{\End}{\mbox{${\rm End}$}}
\newcommand{\Hom}{\mbox{${\rm Hom}$}}
\newcommand{\Ker}{\mbox{${\rm Ker}$}}
\renewcommand{\Im}{\mbox{${\rm Im}$}}

\newcommand{\ta}{\mbox{${\mbox{$\scriptscriptstyle A$}}$}}
\newcommand{\ms}{\mbox{${\mbox{$\scriptscriptstyle M$}}$}}
\newcommand{\ap}{\mbox{$A_{\mbox{$\scriptscriptstyle P$}}$}}
\newcommand{\tx}{\mbox{$\mbox{$\scriptscriptstyle X$}$}}
\newcommand{\tE}{\mbox{$\mbox{$\scriptscriptstyle E$}$}}
\newcommand{\ty}{\mbox{$\mbox{$\scriptscriptstyle Y$}$}}
\newcommand{\kt}{\mbox{$K_{\tx}$}}
\newcommand{\pk}{\mbox{$\pi_{\lambda}^{\scriptscriptstyle K_{X}}$}}
\newcommand{\kk}{\mbox{$K \times_{\scriptscriptstyle K_{X}} (K_{\tx}/T)$}}

\newcommand{\wm}{W^{\tx}}
\newcommand{\wx}{W_{\tx}}
\newcommand{\dw}{\dot{w}}
\newcommand{\fw}{f^{w}_{w_1}}
\newcommand{\lw}{C^{w}_{w_1}}
 
\newcommand{\Dr}{{\rm P}}
\newcommand{\Ulp}{U_{\frak l}^{'}}
\newcommand{\piU}{\pi_{\scriptscriptstyle U}}

\newcommand{\sww}{s^{w}_{w_1}}
\newcommand{\sw}{\Sigma_w}
\newcommand{\swa}{\Sigma_{w_1}}
\newcommand{\swb}{\Sigma_{w_1 w_2}}
\newcommand{\cw}{C_{\dot{w}}}
\newcommand{\cwa}{C_{\dot{w_1}}}

\newcommand{\pp}{\mbox{$\pi_{\mbox{$\scriptscriptstyle P$}}$}}
\newcommand{\pg}{\mbox{$\pi_{\mbox{$\scriptscriptstyle G$}}$}}
\newcommand{\pge}{\mbox{$\pi_{\Delta_{+}}$}}
\newcommand{\iep}{\mbox{$i_{\exp_{\wedge} \pi}$}}
\newcommand{\ienp}{\mbox{$i_{\exp_{\wedge}(- \pi)}$}}

\newcommand{\Xa}{\mbox{$X_{\alpha}$}}
\newcommand{\Ya}{\mbox{$Y_{\alpha}$}}
\newcommand{\pix}{\mbox{$\pi_{\tx,\lambda}$}}
\newcommand{\fix}{\mbox{$\fl_{\tx,\lambda}$}}
\newcommand{\fll}{\mbox{$\fl_{\lambda}$}}

\newcommand{\flds}{\fl_{d, \sigma}}

\newcommand{\bix}{\mbox{$b_{\tx,\lambda}$}}
\newcommand{\dix}{\mbox{$d_{\tx,\lambda}$}}
\newcommand{\pxl}{\mbox{$\partial_{\tx, \lambda}$}}
\newcommand{\tml}{\mbox{$\tilde{m}_{\lambda}$}}
\newcommand{\tse}{\mbox{$T^{*}_{e}(K/T)$}}
\newcommand{\te}{\mbox{$T_{e}(K/T)$}}
\newcommand{\pinf}{\mbox{$\pi_{\infty}$}}
\newcommand{\ponf}{\mbox{$\partial_{\infty}$}}
\newcommand{\ii}{\mbox{$I_{\infty}$}}
\newcommand{\iis}{\mbox{$I_{\infty}^{*}$}}
\newcommand{\il}{\mbox{$I_{\lambda}$}}
\newcommand{\ils}{\mbox{$I_{\lambda}^{*}$}}
\newcommand{\pl}{\mbox{$\pi_{\lambda}$}}
\newcommand{\pal}{\mbox{$\partial_{\lambda}$}}

\newcommand{\hus}{\mbox{$H_{S^1}(K/T, \C)$}}
\newcommand{\hls}{\mbox{$H^{S^1}(K/T, \C)$}}
\newcommand{\cda}{\mbox{$\sigma^{(w_1)}$}}
\newcommand{\cdb}{\mbox{$\sigma^{(w_2)}$}}
\newcommand{\cd}{\mbox{$\sigma^{(w)}$}}
\newcommand{\rank}{{\rm rank}}

\newcommand{\dfl}{\mbox{$d_{{\frak l}_{\lambda}} $}}
\newcommand{\bfl}{\mbox{$b_{{\frak l}_{\lambda}} $}}

\newcommand{\el}{\mbox{$e^{- \lambda}$}}
\newcommand{\rhocheckone}{\mbox{${\check{\rho}}_1$}}
\newcommand{\pkl}{\mbox{$\pi^{{\scriptscriptstyle K_{X}}}_{\lambda}$}}

\newcommand{\dpi}{\mbox{$\delta_{\pi}$}}

\newcommand{\asemi}{\mbox{$\ap \#_{\sigma} A^*$}}
\newcommand{\dsemi}{\mbox{$A \#_{\Delta} A^*$}}

\newcommand{\semi}{\mbox{$\times_{{\frac{1}{2}}}$}}
\newcommand{\fc}{\mbox{${\frak c}$}}
\newcommand{\fu}{\mbox{${\frak u}$}}
\newcommand{\fus}{{\frak u}^*}
\newcommand{\ful}{\fu_{\frak l}}
\newcommand{\flfl}{[\fl, \, \fl]}
\newcommand{\fdl}{\fd_{\frak l}}
\newcommand{\fdlc}{\fd_{\frak l}^{\circ}}

\newcommand{\fd}{\mbox{${\frak d}$}}
\newcommand{\fe}{\mbox{${\frak e}$}}
\newcommand{\fa}{\mbox{${\frak a}$}}
\newcommand{\ft}{\mbox{${\frak t}$}}
\newcommand{\fk}{\mbox{${\frak k}$}}
\newcommand{\fg}{\mbox{${\frak g}$}}
\newcommand{\fq}{\mbox{${\frak q}$}}
\newcommand{\fl}{\mbox{${\frak l}$}}
\newcommand{\fs}{\mbox{${\frak s}$}}
\newcommand{\fsl}{\mbox{${\fs\fl}$}}
\newcommand{\flp}{\mbox{${\frak l}_{p}$}}
\newcommand{\fh}{\mbox{${\frak h}$}}
\newcommand{\fn}{\mbox{${\frak n}$}}
\newcommand{\bfny}{\mbox{${\bar{\frak n}}_{\ty}$}}

\newcommand{\fo}{\mbox{${\frak o}$}}
\newcommand{\fp}{\mbox{${\frak p}$}}
\newcommand{\fr}{\mbox{${\frak r}$}}
\newcommand{\fb}{\mbox{${\frak b}$}}
\newcommand{\fz}{\mbox{${\frak z}$}}
\newcommand{\fm}{\mbox{${\frak m}$}}
\newcommand{\fbp}{\mbox{${\frak b}_{+}$}}
\newcommand{\fbm}{\mbox{${\frak b}_{-}$}}
\newcommand{\fbpm}{\mbox{${\frak b}_{\pm}$}}
\newcommand{\fnp}{\mbox{${\frak n}_{+}$}}
\newcommand{\fnm}{\mbox{${\frak n}_{-}$}}
\newcommand{\fgs}{\mbox{${\frak g}^*$}}
\newcommand{\wg}{\mbox{$\wedge {\frak g}$}}
\newcommand{\wgs}{\mbox{$\wedge {\frak g}^*$}}
\newcommand{\eps}{\mbox{${\epsilon}$}}
\newcommand{\fgdeta}{\mbox{$\fg_{d,\eta}$}}
\newcommand{\fpeta}{\mbox{$\fp_\eta$}}
\newcommand{\fmeta}{\mbox{${\fm}_{\eta}$}}
\newcommand{\dfmeta}{\mbox{${\fm_{\eta}}_{\Delta}$}}
\newcommand{\fnetaplus}{\mbox{${\fn_\eta}_+$}}
\newcommand{\fnetaminus}{\mbox{${\fn_\eta}_-$}}
\newcommand{\fnsigmaplus}{\mbox{${\fn_\sigma}$}}
\newcommand{\fnsigmaminus}{\mbox{${\fn_\sigma}_-$}}
\newcommand{\fnDeltaplus}{\mbox{${\fn_\Delta}_+$}}
\newcommand{\fnDeltaminus}{\mbox{${\fn_\Delta}_-$}}
\newcommand{\fmsDelta}{\mbox{${{\fm}_\Delta}_s$}}
\newcommand{\Geta}{\mbox{$G_\eta$}}
\newcommand{\Peta}{\mbox{$P_\eta$}}
\newcommand{\Gmin}{\mbox{${G^{-\tau_d}}$}}

\newcommand{\fsk}{\mbox{${\fs}_{\frak k}$}}
\newcommand{\As}{\mbox{${\rm Aut}_{\frak g}$}}
\newcommand{\Ints}{\mbox{${\rm Int}_{\frak g}$}}
\newcommand{\Outs}{\mbox{${\rm Aut}_{D({\frak g})}$}}
\newcommand{\daut}{\mbox{$d$}}
\newcommand{\dgamma}{\mbox{${\gamma}_d$}}
\newcommand{\dtheta}{\mbox{$\tau_d$}}
\newcommand{\dDelta}{\mbox{${\Delta,d}$}}
\newcommand{\Zdaut}{\mbox{${Z_d}$}}
\newcommand{\ZdR}{Z_{d, {\Bbb R}}}
\newcommand{\Zdeta}{Z_{d, \eta}}
\newcommand{\Zdetao}{Z_{d, \eta_1}}

\newcommand{\fpS}{\mbox{$\fp_{\scriptscriptstyle S}$}}
\newcommand{\fmS}{\mbox{${\fm}_{\scriptscriptstyle S}$}}
\newcommand{\fnS}{\mbox{${\fn}_{\scriptscriptstyle S}$}}
\newcommand{\fnSone}{\mbox{${\fn}_{\scriptscriptstyle S_1}$}}
\newcommand{\fzS}{\mbox{${\fz}_{\scriptscriptstyle S}$}}
\newcommand{\fzSprime}{\mbox{${\fz}_{{\scriptscriptstyle S}^\prime}$}}
\newcommand{\fzSR}{\mbox{${
{\fz}_{\scriptscriptstyle S}}_{\scriptscriptstyle {\Bbb R}}$}}
\newcommand{\fzSC}{\mbox{${
{\fz}_{\scriptscriptstyle S}}_{\scriptscriptstyle {\Bbb C}}$}}
\newcommand{\fmSss}{\mbox{${\fm}_{\scriptscriptstyle S, 1}$}}
\newcommand{\fmSigss}{\mbox{${\fm}_{\sigma, 1}$}}
\newcommand{\fgreal}{\mbox{${\fg}_0$}}
\newcommand{\PS}{\mbox{$P_{\scriptscriptstyle S}$}}
\newcommand{\NS}{\mbox{$N_{\scriptscriptstyle S}$}}
\newcommand{\MS}{\mbox{$M_{\scriptscriptstyle S}$}}

\newcommand{\cR}{\cal R}
\newcommand{\Meta}{\mbox{$M_\eta$}}
\newcommand{\Netaplus}{\mbox{${N_\eta}+$}}
\newcommand{\Netaminus}{\mbox{${N_\eta}_-$}}
\newcommand{\Cleta}{\mbox{${{\C}^l}_{\eta}$}}
\newcommand{\Rleta}{\mbox{${{\R}^l}_{\eta}$}}
\newcommand{\Rl}{\mbox{${{\R}^l}$}}
\newcommand{\cstar}{\mbox{${\C}^*$}}
\newcommand{\fmsdeltaeta}{\mbox{$(\fmsDelta)^\eta$}}
\newcommand{\Lagr}{\mbox{${\cal L}$}}
\newcommand{\Lagrfd}{\Lagr(\fd)}
\newcommand{\Gr}{{\rm Gr}(n, \fd)}

\newcommand{\LagrfzS}{\mbox{${\cal L}_{{\frak z}_{\scriptscriptstyle S}}$}}
\newcommand{\LagrfzSC}{\mbox{${\cal L}_{{\frak z}_{{\scriptscriptstyle S},
{\Bbb C}}}$}}
\newcommand{\LagrfzSplus}{\mbox{$
{\cal L}_{{\frak z}_{\scriptscriptstyle S,1}}$}}
\newcommand{\LagrfzSminus}{\mbox{$
{\cal L}_{{\frak z}_{\scriptscriptstyle S,-1}}$}}
\newcommand{\LagrfzSeps}{\mbox{$
{\cal L}_{{\frak z}_{\scriptscriptstyle S,\epsilon}}$}}
\newcommand{\zeLagr}{\mbox{${\Lagr}_0$}}
\newcommand{\kLagr}{\mbox{${\Lagr}_{{\frak k}}$}}

\newcommand{\lag}{\Lagr}
\newcommand{\kag}{\kLagr}
\newcommand{\supp}{{\rm supp}}
\newcommand{\fmst}{\fm_{\sigma}^{\theta_{d,\sigma}}}

\newcommand{\Vt}{\mbox{$V_{\ft}$}}
\newcommand{\Va}{\mbox{$V_{\fa}$}}
\newcommand{\Vcplx}{\mbox{$V_0$}}

\newcommand{\wxl}{\mbox{$x_1 \wedge x_2 \wedge \cdots \wedge x_l$}}
\newcommand{\wxk}{\mbox{$x_1 \wedge x_2 \wedge \cdots \wedge x_k$}}
\newcommand{\wyl}{\mbox{$y_1 \wedge y_2 \wedge \cdots \wedge y_l$}}
\newcommand{\wxkm}{\mbox{$x_1 \wedge x_2 \wedge \cdots \wedge x_{k-1}$}}
\newcommand{\wxik}{\mbox{$\xi_1 \wedge \xi_2 \wedge \cdots \wedge \xi_k$}}
\newcommand{\wxikm}{\mbox{$\xi_1 \wedge \cdots \wedge \xi_{k-1}$}}
\newcommand{\wetal}{\mbox{$\eta_1 \wedge \eta_2 \wedge \cdots \wedge \eta_l$}}

\newcommand{\winv}{\mbox{$(\wedge \fg_{1}^{\perp})^{\fg_1}$}}
\newcommand{\wetak}{\mbox{$\eta_1 \wedge \cdots \wedge \eta_k$}}
\newcommand{\gonep}{\mbox{$\fg_{1}^{\perp}$}}
\newcommand{\wonep}{\mbox{$\wedge \fg_{1}^{\perp}$}}

\newcommand{\cala}{\mbox{${\cal A}$}}
\newcommand{\calv}{\mbox{${\cal V}$}}
\newcommand{\pdp}{\mbox{$\partial_{\pi}$}}

\newcommand{\db}{\mbox{$\fd = \fg \bowtie \fgs$}}
\newcommand{\fds}{\mbox{${\scriptscriptstyle {\frak d}}$}}

\newcommand{\Gs}{\mbox{$G^*$}}
\newcommand{\pis}{\mbox{$\pi_{\sigma}$}}
\newcommand{\ea}{\mbox{$E_{\alpha}$}}
\newcommand{\eb}{\mbox{$E_{-\alpha}$}}
\newcommand{\Bm}{\mbox{$ {}^B \! M$}}
\newcommand{\kBm}{\mbox{$ {}^B \! M^k$}}
\newcommand{\Bb}{\mbox{$ {}^B \! b$}}
\newcommand{\epe}{\mbox{$\epsilon$}}
\newcommand{\eot}{\mbox{${\epe \over 2}$}}

\newcommand{\cfg}{\mbox{$C(\fg \oplus \fgs)$}}
\newcommand{\ps}{\mbox{$\pi^{\#}$}}
\newcommand{\backl}{\mathbin{\vrule width1.5ex height.4pt\vrule height1.5ex}}
 
\newcommand{\bx}{\mbox{${\bar{x}}$}}
\newcommand{\by}{\mbox{${\bar{y}}$}}
\newcommand{\bz}{\mbox{${\bar{z}}$}}
\newcommand{\pgs}{\mbox{${\pi_{\mbox{\tiny G}^{*}}}$}}

\newcommand{\tlp}{\mbox{$\tilde{\pi}$}}
\newcommand{\tp}{\mbox{$\tilde{\pi}$}}
\newcommand{\sn}{\mbox{$s_{\scriptscriptstyle N}$}}
\newcommand{\tn}{\mbox{$t_{\scriptscriptstyle N}$}}
\newcommand{\sm}{\mbox{$s_{\scriptscriptstyle M}$}}
\newcommand{\tm}{\mbox{$t_{\scriptscriptstyle M}$}}
\newcommand{\en}{\mbox{$\epsilon_{\scriptscriptstyle N}$}}
\newcommand{\mem}{\mbox{$\epsilon_{\scriptscriptstyle M}$}}

\newcommand{\C}{\mbox{${\Bbb C}$}}
\newcommand{\Z}{\mbox{${\bf Z}$}}

\newcommand{\R}{\mbox{${\Bbb R}$}}
\renewcommand{\a}{\mbox{$\alpha$}}

\newcommand{\flhs}{\fl_{H, \sigma}}
\newcommand{\Grg}{{\rm Gr}(n, \fg)}

 \begin{abstract}
We study Lagrangian subalgebras of a semisimple Lie algebra with
respect to the imaginary part of the Killing form.
We show that
the variety $\Lagr$ of Lagrangian
subalgebras carries a natural 
Poisson structure $\Pi$. We determine the irreducible components of $\Lagr$,
and we show that each irreducible component 
is a smooth fiber bundle over
a generalized flag variety, and  that
the fiber is the product of the real points
of a De Concini-Procesi compactification and a compact homogeneous
space. 
We study some properties of the Poisson structure $\Pi$
and show that it
contains many interesting Poisson submanifolds.
 \end{abstract}

\section{Introduction}
\label{sec_intro}

Let  $\fg$ be a complex semi-simple Lie algebra and let
${\rm Im} \ll\, , \, \gg$ be the  imaginary part of the Killing form
$\ll\, , \, \gg$ of $\fg$.
We will say that a real
subalgebra $\fl$ of $\fg$ is {\it Lagrangian} if
$\dim_{\Bbb R} \fl = \dim_{\Bbb C} \fg$ and if
${\rm Im} \ll x, \,  y\gg = 0$ for all $x, y \in \fl$. 

In this paper, we study the geometry 
of the variety $\Lagr$ of Lagrangian
subalgebras of $\fg$ and show that $\Lagr$ carries a natural 
Poisson structure $\Pi$. We show that each irreducible component of $\Lagr$
is smooth and is a fiber bundle over
a generalized flag variety, and the fiber is the product of the real points
of a De Concini-Procesi compactification and a compact homogeneous
space. 
We study some properties of the Poisson structure $\Pi$
and show that it
contains many interesting Poisson submanifolds.

The Poisson structure $\Pi$ is defined
using the fact that $\fg$, regarded as a real Lie algebra, is the
double of a Lie bialgebra structure on a compact
real form $\fk$ of $\fg$. The construction of 
$\Pi$ works for any Lie bialgebra, and we
present it in the first part of the paper. 
In the second part, we study the specific example of $\Lagr$,
which we regard as the most important example since it is closely
related to interesting problems in Lie theory.

We now explain our motivation and give more details of our results.

Let $(\fu, \fu^*)$ be any Lie bialgebra, let $\fd$ be its double,
and let $\la \, , \, \ra$ be the symmetric scalar product on
$\fd$ given by
\[
\la x + \xi, \, y + \eta\ra \, = \, (x, \, \eta ) \, + \, (y, \, \xi),
\hspace{.2in} x, y \in \fu, \xi, \eta \in \fu^*.
\]
A subalgebra
$\fl$ of $\fd$ is said to be {\it Lagrangian} if $\dim \fl 
= \dim \fu$ and if $\la a, b \ra = 0$ for all
$a, b \in \fl$. Denote by $\Lagrfd$ the set of all
Lagrangian subalgebras of $\fd$. It is a subvariety 
of the Grassmannian of $n$-dimensional subspaces of
$\fd$, where $n = \dim \fu$.
The motivation for studying $\Lagrfd$ comes 
from a
theorem of Drinfeld \cite{dr:homog} on
Poisson homogeneous spaces which we  now recall briefly.
More details are given in Section \ref{sec_drinfi-general}.

Let $(U, \piU)$ be a Poisson Lie group with $(\fu, \fu^*)$
as its tangent Lie bialgebra. Recall that an action of
$U$ on a Poisson manifold $(M, \pi)$ is called Poisson if the 
action map $U \times M \rightarrow M$ is a Poisson map. 
When the action is also transitive, $(M, \pi)$ is called
a $(U, \piU)$-homogeneous Poisson space.
In this case, Drinfeld \cite{dr:homog} associated to each $m \in M$ 
 a Lagrangian subalgebra
$\fl_m$ of $\fd$ and showed that
$\fl_{u \cdot m} = {\rm Ad}_u \fl_m$
for every $u \in U$ and $m \in M$.
 Thus we have a $U$-equivariant map
\begin{equation}
\label{eq_Dr}
P: \, M \lrw \Lagrfd: \, m \Map \fl_m,
\end{equation}
where $U$ acts on $\Lagrfd$ by the Adjoint action.
Drinfeld's theorem says that 
the assignment that assigns to each $(M, \pi)$ the image of
the map $P$ in (\ref{eq_Dr}) gives a one-to-one
correspondence between the set of $U$-equivariant isomorphism
classes of $(U, \piU)$-homogeneous Poisson spaces with 
connected stabilizer subgroups and the set of $U$-orbits in 
a certain subset $\Lagrfd_C$ of $\Lagrfd$ (see Section 
\ref{sec_drinfi-general} for more details).

We prove the following theorem.

\begin{thm}
\label{thm_intro-Pi}

1) There is a Poisson structure $\Pi$ on $\Lagrfd$ with respect to which 
the Adjoint action of $U$ on $\Lagrfd$ is Poisson;

2) Each $U$-orbit ${\cal O}$ in $\Lagrfd$ is a Poisson submanifold
and consequently a $(U, \piU)$-homogeneous Poisson space;

3)  For any $(U, \piU)$-homogeneous Poisson space $(M, \pi)$, the
map $P$ in (\ref{eq_Dr}) 
is a Poisson map onto the $U$-orbit of $\fl_m$ for any $m \in M$.
\end{thm}

We introduce the notation of {\it model points} in $\Lagrfd$.
For a homogeneous Poisson space $(M, \pi)$, let $\fl = P(m)$ for some
$m\in M$. We show $\fl$ is a model point if and only if
 the map $P:M\to {\cal O}_{\frak l}=U
\cdot \fl$
is a local diffeomorphism (and thus a covering map). 
When this happens, we regard 
$({\cal O}_{\frak l},\Pi)$
as a model for the Poisson space $(M, \pi)$. 

The second part of  the paper is concerned with the variety 
$\Lagr$ of
Lagrangian subalgebras of a semi-simple Lie algebra $\fg$ with
respect to the imaginary part of its Killing form.
Let $G$ be
 the adjoint group of
$\fg$. Based on the Karolinsky classification of Lagrangian 
subalgebras of $\fg$ in  \cite{karo:homog-compact}, we prove 

\begin{thm}
\label{thm_introsmooth}
{\it 
The irreducible components of $\Lagr$ are smooth. Each irreducible
component fibers over a generalized flag variety, and its fiber is the
product of a homogeneous space and the space of
 real points of a De Concini-Procesi
compactification of the semisimple part of a Levi subgroup of $G$.
}
\end{thm}

For example, when $\fg = \fsl (2,\C)$, there are two irreducible
components: the first component is the $SL(2,\C)$-orbit through
$\fa + \fn$ and is isomorphic to $\C P^1$ (here $\fa$ consists
of diagonal real trace zero matrices and $\fn$ strictly upper triangular
matrices), and the second component contains the $SL(2,\C)$-orbits
through $\fs\fu (2)$ and $\fsl (2,\R)$ as open orbits, and 
the $SL(2,\C)$-orbit through $i\fa + \fn$ as the unique closed
orbit. The second component may be identified as $\R P^3$. 

Let $\fk$ be a compact real form of $\fg$ and $K \subset G$  the connected
subgroup with Lie algebra $\fk$. 
Then there is a natural Poisson structure $\pi_K$ on $K$
making $(K, \pi_K)$ into a Poisson Lie group such that the 
double of its tangent Lie bialgebra is $\fg$.
By Theorem \ref{thm_intro-Pi}, each $K$-orbit in $\Lagr$ is a
$(K, \pi_K)$-homogeneous Poisson space, and every 
$(K, \pi_K)$-homogeneous Poisson space maps onto a $K$-orbit in $\Lagr$
by a Poisson map. In particular, we show that every point in
the (unique) irreducible component $\Lagr_0$ of $\Lagr$ that
contains $\fk$ is a model point. Consequently, a number of
interesting $(K, \pi_K)$-homogeneous Poisson
spaces are contained in $\Lagr_0$ (possibly up to 
covering maps) as Poisson submanifolds. Among these are
{\it all} $(K, \pi_K)$-homogeneous Poisson structures on 
any $K/K_1$, where $K_1$ is a closed subgroup of $K$
containing a maximal torus of $K$. For example,
$K/K_1$ could be any flag variety $G/Q \cong K/K \cap Q$, where $Q$ is
a parabolic subgroup of $G$.
We remark that it is shown in \cite{lu:cdyb} that all
$(K, \pi_K)$-homogeneous Poisson structures on $K/T$, where $T$
is a maximal torus in $K$, can be 
obtained from solutions to the {\it Classical Dynamical
Yang-Baxter Equation} \cite{e-v:cdyb}. Some Poisson
geometrical properties of such Poisson structures are
also studied in \cite{lu:cdyb}.

We are motivated to study $(K, \pi_K)$-homogeneous Poisson structures 
because of their connections to Lie theory. One remarkable example is
the
so-called Bruhat Poisson structure
$\pi_{\infty}$ \cite{lu-we:poi}
on $K/T$. It corresponds
to the Lagrangian subalgebra 
$\ft+\fn$ of $\fg$, where $\fg=\fk+\fa+\fn$
is an Iwasawa decomposition of $\fg$, and $\ft=i\fa$ is the Lie algebra
of $T$.  The name Bruhat Poisson structure
 comes from the fact that its symplectic leaves are exactly
the Bruhat cells for a Bruhat decomposition of $K/T$ 
\cite{lu-we:poi}; its
Poisson cohomology is isomorphic to a direct sum of 
$\fn$-cohomology groups with coefficients in certain
principal representations of $G$
\cite{lu:homog}; its $K$-invariant {\it Poisson harmonic forms}
are exactly the harmonic forms introduced and studied by Kostant
in \cite{ko:63}. This last fact is proved in \cite{e-l:harm},
where we also use $\pinf$ to construct $S^1$-equivariantly
closed forms on $K/T$ and use them to reinterpret
the Kostant-Kumar approach to the Schubert calculus on $K/T$
\cite{k-k:integral}. One key fact used in \cite{e-l:harm} is that
the Poisson
structure $\pi_{\infty}$ is the limit of a family $\pi_t$,
$t \in (0, +\infty)$, of 
$(K,\pi_K)$-homogeneous symplectic structures 
on $K/T$.
The family $\pi_t$ corresponds to a continuous
curve in $\Lagr$.
Thus, we regard $\Lagr$ as a natural setting for
deformation problems for Poisson homogeneous spaces, and for this
reason it is desirable to study its geometry.

The paper is organized as follows.

We start our discussion in Section \ref{sec_general}
with an arbitrary Poisson Lie group $(U, \piU)$, its
tangent  Lie bialgebra $(\fu, \fus)$, 
and the variety $\Lagr(\fd)$ of Lagrangian subalgebras
of its double $\fd = \fu \bowtie \fus$. We first review Drinfeld's
theorem on $(U, \piU)$-homogeneous spaces. 
We then give the construction of the Poisson structure
$\Pi$ on $\Lagr(\fd)$ and establish the
properties listed in Theorem \ref{thm_intro-Pi}.

The rest of the paper is devoted to  
the Poisson
Lie group $(K, \pi_K)$. 
In \ref{sec_kara-classification},
 we review Karolinsky's classification of Lagrangian subalgebras,
and use it to decompose
$\Lagr$ into a finite  disjoint union of submanifolds $\Lagr(S,\eps,d)$.
The study of the closure $\overline{\Lagr(S,\eps,d)}$
is reduced to studying the closure of the variety of real forms
of a semisimple Lie algebra. After some preliminary results in
Section \ref{sec_sig}, we identify the closure with the real points
of a De Concini-Procesi compactification in Section \ref{sec_dp}.
In Section \ref{sec_lsedclos}, we apply our results
to determine the irreducible components of $\Lagr$ and show they
are smooth. We also study the set of model points in $\Lagr$ and show
that every Lie algebra in the irreducible component ${\cal L}_0$
containing $\fk$
is a model point. 
Finally, in Section \ref{sec_poi-on-lagr}, we study some 
properties of the Poisson structrure $\Pi$. In particular, we study the 
$K$-orbits in the irreducible component $\Lagr_0$ and the 
$(K, \pi_K)$-homogeneous Poisson spaces arising from them.

We would like to thank Eugene Karolinsky and Hermann Flaschka
for useful conversations,
and the Banach Center for its hospitality when some of these results
were found. In addition, the first author would like to thank
Northwestern University and the University of Chicago  and
the second author the Hong Kong University of Science and Technology
for their hospitality during the preparation of the paper.

\section{Generalities on Lie bialgebras}
\label{sec_general}

\subsection{Drinfeld's theorem}
\label{sec_drinfi-general}

In this section, we review Drinfeld's theorem on 
homogeneous spaces of Poisson Lie groups in  \cite{dr:homog}.
Details on Poisson Lie groups can be found in
\cite{lu-we:poi} and \cite{k-s:quantum} and the references cited
in \cite{k-s:quantum}.

Let $(U, \piU)$ be a Poisson Lie group with tangent
Lie bialgebra $(\fu, \fus)$, where $\fu$ is the Lie algebra of
$U$ and $\fus$ its dual space equipped with a Lie algebra
structure coming from the linearization of $\piU$ at the 
identity element of $U$. We will use letters $x, y, x_1, y_1, 
\cdots$ to denote elements in $\fu$ and $\xi, \eta, \xi_1, \eta_1, 
\cdots$ for elements in $\fus$. The pairing between elements
in $\fu$ and in $\fus$ will be denoted by $( \, , \, )$.

Let $\la \, , \, \ra$ be the symmetric non-degenerate scalar product
on the direct sum vector space $\fu \oplus \fus$ defined by
\begin{equation}
\label{eq_scalar}
\la x_1 + \xi_1, \, x_2 + \xi_2 \ra \, = \, (x_1, \xi_2) + (x_2, \xi_1).
\end{equation}
Then there is a unique Lie bracket on the 
$\fu \oplus \fus$ such that $\la \, , \, \ra$ is ad-invariant
and that both $\fu$ and $\fus$ are its Lie subalgebras with 
respect to the
natural inclusions.  The vector space $\fu \oplus \fus$ together
with this Lie bracket is called
the double Lie algebra of $(\fu, \fus)$ and we will denote it by 
$\fd = \fu \bowtie \fus$.
Note that $U$ acts on $\fd$ by the Adjoint action (by first mapping
$U$ to the adjoint group of $\fd$).

\begin{exam}
\label{exam_lagrkan}
{\em
Let $\fu = \fk$ be a compact semi-simple Lie algebra. Let
$\fg={\fk}_{\Bbb C}$ be the complexification  of $\fk$ with an
Iwasawa decomposition $\fg=\fk+\fa+\fn$. Let $\la \, , \, \ra$
be  twice the imaginary part of the Killing form of $\fg$. 
Then the pairing
between $\fk$ and $\fa + \fn$ via $\la \, , \, \ra$ gives an
identification of $\fk^*$ and $\fa + \fn$, and $(\fk, \fa + \fn)$
becomes a Lie bialgebra whose double is $\fg$. If $K$ is any group
with Lie algebra $\fk$, then there is a Poisson structure $\pi_K$
on $K$ making $(K, \pi_K)$ into a Poisson Lie group whose
tangent Lie bialgebra is $(\fk, \fa + \fn)$.
This will be our most important example.}
\end{exam}

\begin{dfn}
\label{dfn_lagr-general}
{\em Let $n = \dim \fu$. 
 A Lie subalgebra $\fl$ of $\fd$ is called {\it Lagrangian}
if $\la a, \, b \ra = 0$ for all $a, b \in \fl$ and if
$\dim \fl = n$. The set of all Lagrangian subalgebras of $\fd$
will be denoted by $\Lagrfd$.
}
\end{dfn}
 
Both $\fu$ and $\fu^*$ are Lagrangian. If $D$ is the adjoint 
group of $\fd$, then $D$ acts on the set of Lagrangian subalgebras.
In Example \ref{exam_lagrkan}, any real form of $\fg$
is a Lagrangian subalgebra, as is $\ft+\fn$, where $\ft = i \fa$ is the
centralizer of $\fa$ in $\fk$.

\bigskip
Let $(M, \pi)$ be a $(U, \piU)$-homogeneous Poisson space. 
Recall \cite{dr:homog} that this means that
$U$ acts on $M$ transitively and that the action
map $U \times M \rightarrow M$ is a Poisson map, where
$U \times M$ is equipped with the direct product Poisson
structure $\piU \oplus \pi$.
Let $m \in M$.
Then being $(U, \piU)$-homogeneous, the
Poisson structure $\pi$ on $M$ must satisfy
\begin{equation}
\label{eq_poi-action}
\pi(um) \, = \, u_* \pi(m) \, + \, m_* \piU(u), \hspace{.2in}
\forall u \in U, \, m \in M.
\end{equation}
Here $u_*$ and $m_*$ are respectively
the differentials of the
maps $M \rightarrow M: m_1 \mapsto um_1$ and
$U \rightarrow M: u_1 \mapsto u_1 m$. Thus, $\pi$ is
totally determined by its value $ \pi(m) \in \wedge^2
(T_m M)$ at $m$. Let
$U_m \subset U$ be the stabilizer subgroup of $U$ at $m$
with Lie algebra $\fu_m$. 
 Identify
$T_mM \cong \fu /\fu_m$ so that $\pi(m) \in \wedge^2 (\fu /\fu_m)$.
Let $\fl_m$ be the subspace of
$\fd$ defined by
\begin{equation}
\label{eq_flm}
\fl_m \, = \, \{x+\xi: x \in \fu,  \, \xi \in \fus,
\xi|_{{\frak u}_m} = 0, \xi \backl \pi(m) = x + \fu_m\}.
\end{equation}

\begin{thm}[Drinfeld \cite{dr:homog}]
\label{thm_drinfi-homog}
1) $\fl_m$ is a Lagrangian subalgebra of $\fd$ for all $m \in M$;

2) For all $m \in M$ and $u \in U$,
\begin{eqnarray}
\label{eq_drinfi-1}
& & \fl_m \cap \fu \,  = \,  \fu_m \\
\label{eq_drinfi-2}
& & {\rm Ad}_u \fl_m \, = \, \fl_{um}, \hspace{.2in}
\forall u \in U.
\end{eqnarray}

3) Let $M$ be a $U$-homogeneous space. A $(U, \piU)$-homogeneous
Poisson structure $\pi$ on $M$ is equivalent to a 
$U$-equivariant map $\Dr: M \rightarrow \Lagrfd: m \mapsto \fl_m$ such
that (\ref{eq_drinfi-1}) holds
for all $m \in M$.
\end{thm}

\begin{dfn}
\label{dfn_drinfi-lag}
{\em
We will call $\fl_m$ the {\it Lagrangian subalgebra of $\fd$
associated to $(M, \pi)$ at the point $m$}.
The map $\Dr: M \rightarrow \Lagrfd$ will be called the {\it Drinfeld map}.
}
\end{dfn}

\begin{dfn}
\label{dfn_property-I}
{\em
Given a $U$-homogeneous space $M$, we say that a
$U$-equivariant map $M \rightarrow \Lagrfd: m \mapsto \fl_m$
has 
{\it Property I} (I for intersection) if 
(\ref{eq_drinfi-1}) is satisfied  for all $m \in M$.
}
\end{dfn}

Thus 3) of Theorem \ref{thm_drinfi-homog} can be rephrased as
follows: given a $U$-homogeneous space $M$, a 
$(U, \piU)$-homogeneous Poisson structure on $M$ is equivalent to 
a $U$-equivariant map $M \rightarrow \Lagrfd$  with
Property I.

\begin{rem}
\label{rem_bundle-to-poi}
{\em
We explain how a $U$-equivariant map $M \rightarrow \Lagrfd$ 
having Property I gives
a $(U, \piU)$-homogeneous Poisson structure on $M$: pick any $m
\in M$. Because $\fl_m \subset \fd$ is maximal isotropic
(this means that $\dim \fl_m = n$ and that $\la a, \, b\ra = 0
$ for all $a, b \in \fl_m$) and because of (\ref{eq_drinfi-1}),
an easy linear algebra argument (see also
Lemma \ref{lem_linear-algebra-fact}) shows that 
there is a unique element $\pi(m) \in \wedge^2 (\fu / \fu_m)$
such that (\ref{eq_flm}) holds. Define a bivector field
$\pi$ on $M$ by (\ref{eq_poi-action}). This is well defined because
of (\ref{eq_drinfi-2}). This $\pi$ is Poisson  because
$\fl_m$ is Lagrangian.  It is
$(U, \piU)$-homogeneous because (\ref{eq_poi-action}) holds by definition.
}
\end{rem}

We now state some consequences of Theorem \ref{thm_drinfi-homog}.

\begin{dfn}
\label{dfn_C}
{\em
A Lagrangian subalgebra of $\fd$ is said to have {\it Property C}
(C for closed) if the connected subgroup  $\Ulp$ of $U$ with Lie
algebra $\fl \cap \fu$ is closed in $U$.
}
\end{dfn}

Note that any $\fl_m$ in the image of the Drinfeld map
for any $(M, \pi)$ has Property $C$, because the connected subgroup of
$U$ with Lie algebra $\fl_m \cap \fu$ is the identity
connected component
of the stabilizer subgroup of $U$ at $m$, so it is closed in $U$. 
Conversely, if $\fl \in \Lagrfd$
has Property C, we  have the $U$-homogeneous space $U/U_{\frak l}^{'}$
and the $U$-equivariant map 
\[
U/U_{\frak l}^{'} \lrw \Lagrfd: \, uU_{\frak l}^{'} \Map {\rm Ad}_u \fl.
\]
It has Property I. 
More generally, suppose that
$U_1$ is any closed subgroup of $U$ having the properties

A) the Lie algebra of $U_1$ is $\fl \cap \fu$;

B) $U_1$ normalizes
$\fl$,

Then we have the $U$-equivariant map
\[
U/U_1 \lrw \Lagrfd: \, uU_1 \Map {\rm Ad}_u \fl.
\]
It has Property I.
Thus, by Theorem \ref{thm_drinfi-homog},
we have

\begin{cor}
\label{cor_Uone}
Suppose that $\fl \in \Lagrfd$ has Property C. Then for any
closed subgroup $U_1$ of $U$ having Properties A) and B),
 there is a $(U, \piU)$-homogeneous
Poisson structure on $U/U_1$ whose Drinfeld map is given by
\[
\Dr: \, U /U_1 \, \lrw \, \Lagrfd: \, uU_1 \Map 
{\rm Ad}_u \fl.
\]
\end{cor}

\begin{dfn}
\label{dfn_standard-C}
{\em 
For a Lagrangian subalgebra $\fl$ of $\fd$ with Property C and
any subgroup $U_1$ of $U$ with the above Properties A) and B),
we say that the 
Poisson manifold $(U/U_1, \pi)$ described in Corollary \ref{cor_Uone}
is determined by $\fl$. 
}
\end{dfn}

Denote by $\Lagrfd_C$ the set of all points in $\Lagrfd$ with
Property C. It is clearly invariant under the Adjoint action of $U$.
For every  $(U, \piU)$-homogeneous Poisson space $(M, \pi)$, the image
of the Drinfeld map $M \rightarrow \Lagrfd$ is a $U$-orbit in 
$\Lagrfd_C$. 

\begin{cor} [Drinfeld \cite{dr:homog}]
\label{cor_one-to-one}
The map that assigns to each $(M, \pi)$
  the image of its Drinfeld map gives a one-to-one
correspondence between $U$-equivariant isomorphism classes of
 $(U, \piU)$-homogeneous Poisson spaces with connected stabilizer subgroups
and the set of $U$-orbits in $\Lagrfd_C$.
\end{cor}


We close this section by an example of a Lagrangian subalgebra  $\fl$
that does not have Property C.

\begin{exam}
\label{exam_notclosed}
{\em
\cite{karo:homog-compact} Consider the Lie bialgebra
 $(\fk, \fa + \fn)$ in Example 
\ref{exam_lagrkan}. Let $U=K$ be a compact connected
Lie group with Lie algebra $\fk$ and let
 $T$ be the maximal torus of $K$ with Lie algebra $i\fa$.
 Choose a topological generator $t$ of $T$
and let $t=\exp(X), X\in \ft$. 
Let $\fl=\R\cdot X + (\fa \cap (\R\cdot X)^{\perp})
+ \fn$, where the perpendicular is computed relative to the Killing form.
Then $\fl$ is Lagrangian, but if ${\rm rank}(T) > 1$ then $\fl\cap \fk$ is not
the Lie algebra of a closed subgroup of $K$, so $\fl$ does not
have Property C. 
}
\end{exam}

\subsection{A ``Poisson structure" on $\Lagrfd$}
\label{sec_Lagr-general}
 
Let $(U, \piU)$ be a Poisson Lie group and let $(\fu, \fus)$ be its
tangent Lie bialgebra. Let $\fd = \fu \bowtie \fus$ be its
double Lie algebra equipped with the symmetric scalar product
$\la \, , \, \ra$ given by (\ref{eq_scalar}).
Recall that $\Lagr(\fd)$ is the set of Lagrangian subalgebras of
$\fd$ with respect to $\la \, , \, \ra$.

\begin{nota}
\label{nota_lagr-general}
{\em
We will 
use ${\rm Gr}(n, \fd)$
to denote the Grassmannian of $n$-dimensional
subspaces of $\fd$. Since the condition of being closed under
Lie bracket and the condition of being Lagrangian are polynomial
conditions, $\Lagrfd \subset \Gr$ is an algebraic
subset. 
}
\end{nota}

The group $U$ acts on $\Gr$ by the Adjoint action and it leaves 
$\Lagrfd$ invariant. 
Although $\Lagrfd$ may be singular, all the $U$-orbits in
$\Lagrfd$ are smooth.

In this section, we will show that there is a smooth bi-vector
field $\Pi$ on $\Gr$ with the property
\[
[\Pi, \, \Pi](\fl) \, = \, 0
\]
for every $\fl \in \Lagrfd$, where $[\Pi, \, \Pi]$ is the 
Schouten bracket of $\Pi$ with itself. Moreover,
we show that  $\Pi$ is
tangent to every $U$-orbit $\cO$ in $\Lagrfd$, so 
$(\cO, \Pi)$ is a Poisson manifold. In fact, each $(\cO, \Pi)$
is a $(U, \piU)$-homogeneous Poisson space.
If $(M,\pi)$ is a $(U, \piU)$-homogeneous Poisson space,
we show that the Drinfeld map $\Dr: M \rightarrow {\cal O}$ 
is a Poisson map, where ${\cal O}$ is the 
 $U$-orbit of $\fl_m$ for any $m \in M$.

\begin{nota}
\label{nota_sharp}
{\em
We identify $\fd^* \cong \fus \oplus \fu$
in the obvious way. Denote by $\# :\fd^* \to \fd$ the isomorphism 
induced by the nondegenerate pairing $\la \, , \, \ra$ on $\fd$. It is 
given by
\begin{equation}
\label{eq_fds-identify}
\#: \, \fd^* \lrw \fd: \, \#(\xi + x)=x+\xi.
\end{equation}
For $V\subset \fd$, we let 
\[
V^\circ=\{ f\in \fd^*: \,  f|_V = 0 \}.
\]
}
\end{nota}

To define the bi-vector field $\Pi$ on $\Gr$, we consider
 the 
element 
$R \in \wedge^2 \fd$ 
defined by 
\[
 R( \xi_1 + x_1, \, \xi_2 + x_2) \, = \, (\xi_2, \, x_1) \, - \, 
(\xi_1, \, x_2), \hspace{.2in} \forall x_1, x_2 \in \fu, \, 
\xi_1, \xi_2 \in \fus.
\]
The element $R$ is an example of a classical $r$-matrix on $\fd$
\cite{k-s:quantum}. In particular, the Schouten bracket 
$[R, R] \in \wedge^3 \fd$ of $R$ with itself
is ad-invariant and is given by
\beqa
[R,R](f_1,f_2,f_3)=2<\# f_1, \, [\# f_2,\# f_3]>
\eeqa
for $f_i\in \fd^*$. Denote by $\chi^k(\Gr)$ the space of $k$-vector fields on
$\Gr$ (i.e., the space of smooth sections of
the $k$-th exterior power of the tangent bundle of $\Gr$).
The action by the adjoint group $D$ of $\fd$ on $\Gr$ gives a Lie
algebra anti-homomorphism
\[
\kappa: \, \fd \lrw \chi^1(\Gr)
\]
whose multi-linear extension from $\wedge^k \fd$ to $\chi^k(\Gr)$,
for any integer $k \geq 1$,
will also be denoted by $\kappa$.

Define the bi-vector field $\Pi$ on $\Gr$ by
\[
\Pi \, = \, {1 \over 2} \kappa (R).
\]

\begin{thm}
\label{thm_Pi-0}
For every Lagrangian subalgebra $\fl$ of $\fd$ regarded
as a point in $\Gr$, we have
\[
[\Pi, \, \Pi](\fl) \, = \, 0,
\]
where $[\Pi, \Pi]$ is the Schouten bracket of $\Pi$ with itself.
\end{thm}

\noindent
{\bf Proof.}
Since $\Pi = {1 \over 2} \kappa(R)$ and since $\kappa$
is a Lie algebra anti-homomorphism, we have
\[
[\Pi, \, \Pi] \, = \, - {1 \over 4} \kappa ([R, \, R]).
\]
Let $D_{\frak l}$ be the stabilizer subgroup of $D$ at $\fl$
for the Adjoint action, and let $\fd_{\frak l}$ be its Lie
algebra. Since $\Pi$ is tangent to the $D$-orbit $D \cdot \fl$
in $\Gr$, we only need to show that $[\Pi,\Pi]=0$ when evaluated on
a triple $(\alpha_1,\alpha_2,\alpha_3)$ of covectors in
$T_{\frak l}^*(D\cdot \fl)$. 
The map
\[
\kappa: \, \fd \lrw T_{\frak l}(D \cdot \fl)
\]
gives an identification
\[
\kappa^*: \, T_{\frak l}^{*} (D \cdot \fl) \lrw \fd_{\frak l}^{\circ},
\]
Thus, it suffices to show 
\[
[R,R](f_1,f_2,f_3)=0
\]
for $f_i \in \fd_{\frak l}^{\circ}, i = 1, 2, 3$. Since $\fl \subset
\fd_{\frak l}$,  we have 
$\#(\fd_{\frak l}^{\circ})\subset \#({\fl}^\circ)=\fl$.
It follows that
\[
[R,R](f_1,f_2,f_3)=2<\# f_1,[\# f_2,\# f_3]>=0
\]
because $\fl$ is a Lagrangian subalgebra.

\qed

\begin{cor}
\label{cor_D-orbits}
For every  $\fl \in \Lagrfd \subset \Gr$,
the bivector field $\Pi$ defines
a Poisson structure on the $D$-orbit $D \cdot \fl$ in
$\Gr$.
\end{cor}

Since $[R, R] \in \wedge^3 \fd$ is ad-invariant, the following bivector
field $\pi_-$ on $D$ is Poisson:
\[
\pi_- (d) \, = \, {1 \over 2}( r_d R \, - \, l_d R), \hspace{.2in}
d \in D,
\]
where $r_d$ and $l_d$ are respectively the differentials of the 
right and left translations on $D$ defined by $d$. Moreover, 
$(D, \pi_-)$ is a Poisson Lie group and $(U, \piU)$ is a
Poisson subgroup of $(D, \pi_-)$ (see \cite{lu:thesis}).

\begin{prop}
\label{prop_D-also}
For every  $\fl \in \Lagrfd$,
the Poisson manifold $(D \cdot \fl, \, \Pi)$ is
 $(D, \pi_-)$-homogeneous.
\end{prop}

\noindent
{\bf Proof.}
Let again  $D_{\frak l}$ be the stabilizer subgroup of $\fl$ in $D$. Then
$D \cdot \fl \cong D / D_{\frak l}$. Consider the bivector field
$\Pi_1$ on $D$ defined by 
\[
\Pi_1 (d) \, = \, {1 \over 2} r_d R, \hspace{.2in} d \in D.
\]
Then $\Pi = p_* \Pi_1$, where $p: D \rightarrow D / D_{\frak l}$
is the natural projection and
$p_*$ its differential. It is easy to check that for any
$d_1, d_2 \in D$, we have
\[
\Pi_1(d_1d_2) \, = \, l_{d_1} \Pi_1(d_2) \, + \, r_{d_2} \pi_-(d_1).
\]
It follows that $(D \cdot \fl, \, \Pi)$ is
a $(D, \pi_-)$-homogeneous Poisson space.
\qed

Consider now the $U$-orbits in $\Lagrfd$ through a point
$\fl \in \Lagrfd$.
We have

\begin{thm}
\label{thm_Pi-U-orbits}
At any $\fl \in \Lagrfd$,
the bi-vector field $\Pi$ on $\Gr$
is tangent to the $U$-orbit through $\fl$, so that
$(U \cdot \fl, \, \Pi)$ is a Poisson submanifold of $(D \cdot \fl, \, \Pi)$.
\end{thm}

\noindent
{\bf Proof.}
Regard $\Pi$ as a bivector field on the $D$-orbit 
$D \cdot \fl$, so $\Pi(\fl) \in \wedge^2 T_{\frak l}(D \cdot \fl)$.
Let $\Pi(\fl)^{\#}$ be the linear map
\beqa
\Pi(\fl)^{\#}:  & &  T_{\frak l}^{*} (D \cdot \fl) \lrw 
T_{\frak l}(D \cdot \fl): \\  
& & \Pi(\fl)^{\#}(\alpha) (\beta) = \Pi(\fl)(\alpha, \beta), \hspace{.2in}
\alpha, \beta \in T_{\frak l}^{*} (D \cdot \fl).
\eeqa
It is enough to show that the image of $\Pi(\fl)^{\#}$ is tangent to the 
$U$-orbit through $\fl$.

By the identification,  $T_{\frak l}^{*} (D \cdot \fl) \rightarrow
\fd_{\frak l}^{\circ}$,
it is enough to show that 
\[
\kappa((\xi + x) \backl R) \, \in \,T_{\frak l}(U \cdot \fl), \hspace{.2in}
\forall \xi + x \in \fd_{\frak l}^{\circ},
\]
where $(\xi + x) \backl R \in \fd$ is defined by
\[
((\xi + x) \backl R)(\eta + y) \, = \, R(\xi + x, \, \eta + y),
\hspace{.2in} \forall \eta + y \in \fd^*.
\]
We compute explicitly. It follows from the definition of $R$ that
\[
R \, = \,   \sum_{i=1}^{n} \eta_i \wedge e_i \, \in \, \wedge^2 \fd,
\]
where
$\{e_1, ..., e_n\}$
is a basis for $\fu$ and $\{\eta_1, ..., \eta_n\}$ is
its dual basis for
$\fus$.  It follows that
\[
(\xi + x) \backl R \, = \, \sum_{i=1}^{n} \left(
(x, \eta_i) e_i - (\xi, e_i) \eta_i  \right) \, = \, x - \xi.
\]
Hence
\[
\kappa((\xi + x) \backl R) \, = \, \kappa(x) \, - \, \kappa(\xi).
\]
But since $\xi + x \in \fd_{\frak l}^{\circ}$, we have
$x + \xi \in \fl$, so $\kappa(x + \xi) = 0$. Thus
\[
\kappa((\xi + x) \backl R) \, = \, 2 \kappa(x) \, \in \, 
T_{\frak l}(U \cdot \fl).
\]
\qed
 
 
\begin{cor}
\label{cor_U-homog}
For every $\fl \in \Lagrfd$, the Poisson manifold
$(U \cdot \fl, \, \Pi)$ is a $(U, \piU)$-homogeneous Poisson
space.
\end{cor}

\noindent
{\bf Proof.} This follows from Proposition
\ref{prop_D-also} because $(U, \piU)$ is a Poisson subgroup of
$(D, \pi_-)$ and $(U \cdot \fl, \, \Pi)$ is a Poisson
submanifold of $(D \cdot \fl, \, \Pi)$.
\qed

\begin{rem}
\label{rem_dualgroup}
{\em
Let $U^*$ be the connected  and simply connected group 
with Lie algebra ${\fu}^*$.
Then for any Lagrangian subalgebra $\fl \in \Lagrfd$, the orbit
$U^* \cdot \fl$ is also a Poisson submanifold of $(D \cdot \fl,\Pi)$.
Indeed, the roles of $\fu$ and ${\fu}^*$ are symmetric in the definition
of $D$ and of $\Lagrfd$, but the $R$-matrix for the Lie bialgebra
$(\fu^*, \fu)$ differs from that for $(\fu, \fu^*)$ by a minus sign.
Consequently, if we denote by $\pi_{{\scriptscriptstyle U^*}}$
the Poisosn structure on $U^*$ such that 
$(U^*, \pi_{{\scriptscriptstyle U^*}})$ is the dual Poisson Lie group of 
$(U, \piU)$, then every $U^*$-orbit in $\Lagrfd$ is a
$(U^*, -\pi_{{\scriptscriptstyle U^*}})$-homogeneous Poisson space.
}
\end{rem}

We now  look at the Drinfeld map $\Dr: U \cdot \fl \rightarrow
\Lagrfd$ for the $(U, \piU)$-homogeneous Poisson
space $(U \cdot \fl, \pi)$ (see Definition \ref{dfn_drinfi-lag}).

\begin{thm}
\label{thm_Pi-Lagr-at-fl}
For any $\fl \in \Lagrfd$,
the Lagrangian 
subalgebra of $\fd$ associated to $(U \cdot \fl, \, \Pi)$
at $\fl$ is
\[
T(\fl ) \, = \, \ful \,+ \, (\fu +  \fu_{{\frak l}}^{\perp}) \cap \fl,
\]
where $\ful$ is the normalizer subalgebra of $\fl$ in $\fu$,
and  $\fu_{{\frak l}}^{\perp} = \{\xi \in \fu^*: \, 
\xi|_{\fu_{{\frak l}}} = 0\}.$
\end{thm}

\noindent
{\bf Proof.} 
Denote by
$\fl^{'}$ the Lagrangian subalgebra
associated to $(U \cdot \fl, \, \Pi)$ at $\fl$.  We need to
show that  $\fl^{'} = T(\fl)$.
By definition,
\[
\fl^{'} \, = \, \{x + \xi: \, x \in \fu, \, \xi \in \fu_{\frak l}^{\perp},\,
\xi \backl \Pi(\fl) = x + \ful\}.
\]
Let $\xi \in \fu_{\frak l}^{\perp}$. Since 
the inclusion
\[
(U \cdot \fl, \, \Pi) \lrw (D \cdot \fl, \, \Pi)
\]
is a Poisson map, it suffices to compute $((\kappa^*)^{-1}(\xi + x)) 
\backl \Pi(\fl)$
for any $x\in \fu$ such that $\xi + x \in {\fd}_{\frak l}^0$, where 
$\Pi(\fl)$ is regarded as a bi-vector
at $\fl \in D \cdot \fl$, and $(\kappa^*)^{-1}: T_{{\frak l}}^{*}(D
\cdot \fl) \rightarrow \fd_{{\frak l}}^{\circ}$ is the 
isomorphism induced by $\kappa: \fd \rightarrow 
T_{{\frak l}}(D \cdot \fl)$. 
In the proof of Theorem \ref{thm_Pi-U-orbits}, we showed that
$(\kappa^*)^{-1}(\xi + x) \backl \Pi(\fl) = \kappa(x)$.
As a result, we see that
\beqa
\fl^{'} &  = &  \{x+ \xi:\, \xi + x_1 \in \fd_{\frak l}^{\circ} 
\, \,  {\rm for} \, {\rm some} \, x_1 = x \, {\rm mod} (\ful)\}\\
& = & \ful \, + \,  \#(\fd_{\frak l}^{\circ}).
\eeqa
Now the inclusions
${\fu}_{\frak l}\subset {\fd}_{\frak l}$ and $\fl \subset 
{\fd}_{\frak l}$
induce inclusions $\#({\fd}_{\frak l}^\circ) \subset \fu + 
\fu_{{\frak l}}^{\perp}$ 
and
$\#({\fd}_{\frak l}^\circ) \subset \fl$, so 
$\#({\fd}_{\frak l}^\circ) \subset (\fu + \fu_{{\frak l}}^{\perp})
\cap \fl$. Hence,
\[
\ful \, + \, \#({\fd}_{\frak l}^\circ) \subset
\ful \, + \, 
(\fu + {\fu}_{\frak l}^\perp) \cap \fl = T(\fl).
\]
On the other hand, it is obvious that $T(\fl)$ is isotropic, so
its dimension is at most $n$. Since $\fl^{'}$ has dimension $n$, we 
must have
$\fl^{'} = T(\fl)$.
\qed
 
\begin{rem}
\label{rem_notcont}
{\em
The map $T:\Lagrfd \to \Lagrfd$ is not continuous in general.
For example, consider the Lie bialgebra in Example \ref{exam_lagrkan}
for $\fg = \fsl (3,\C)$.
Choose $H\in \fa$ with the property that both simple roots are
positive on $H$ and consider the curve $\gamma_t = \exp(
{\rm ad}_{tH})(\fsl (3,\R))$
in $\Lagr = \Lagr(\fg)$.
Let $\gamma_{\infty}$ be the limit of $\gamma_t$ as $t\to \infty$
in $\Lagr$.
Clearly, $\gamma_t$ is isomorphic to $\fsl (3,\R)$ for
$t\not= \infty$, and one can show $\gamma_{\infty} = {\fh}^{\tau} + \fn$,
where $\fh = \fa + \ft$ is a Cartan subalgebra of $\fsl (3,\C)$,
and $\tau$ is an anti-linear automorphism such that
$\dim ({\fh}^\tau \cap \ft)=1$. We will show later that
when $\fl$ is a real form of a complex semi-simple Lie algebra,
then $\fl$ is its own normalizer. It follows that $T(\gamma_t)=
\gamma_t$ for all $t < \infty$.
On the other hand, it is easy to check that $T(\gamma_{\infty})=\ft
+ \fn$. It follows that $T$ is not continuous. This example can
be generalized to any real form corresponding to a nontrivial
diagram automorphism 
(see Remark \ref{rem_closedorbit} for a generalization
of this example).
}
\end{rem}

Assume now  that $(M, \pi)$ is an arbitrary $(U, \piU)$-homogeneous 
Poisson space.  Consider the Drinfeld map
\[
\Dr: \, M \lrw \Lagrfd: \, m \Map \fl_m.
\]
By Theorem \ref{thm_drinfi-homog}, $\Dr$ is a submersion of
$M$ onto the $U$-orbit $\cO = U \cdot \fl_m$ in $\Lagrfd$ for any 
$m \in M$.

\begin{thm}
\label{thm_M-to-O}
The Drinfeld map 
\[
\Dr: \, (M, \, \pi) \lrw (\cO, \, \Pi)
\]
is a Poisson map.
\end{thm}

\noindent
{\bf Proof.}
Fix $m \in M$. Let $\fl = \fl_{m}$. Then $\cO = U \cdot \fl$.
Since both $(M, \pi)$ and $(\cO, \Pi)$
are $(U, \piU)$-homogeneous, it is enough to show that
\[
\Dr_* \pi(m) \, = \, \Pi (\fl).
\]
Let $U_{m}$ and $U_{{\frak l}}$ be respectively
the stabilizer subgroup of $U$ at $m$
and the normalizer subgroup of $\fl$
in $U$. Their Lie algebras are respectively
$\fl \cap \fu$ and $\fu_{{\frak l}}$.
Since $\Dr$ is $U$-equivariant, we have $U_{m} \subset U_{{\frak l}}$.
Identify
\[
M \, \cong \, U/U_{m}, \hspace{.2in}
\cO \, \cong \, U/U_{{\frak l}}.
\]
Then the map $\Dr$ becomes
\[
\Dr: \, U/U_{m} \lrw U/U_{{\frak l}}: \, uU_{m} \Map uU_{{\frak l}},
\hspace{.2in}
\]
and we have 
\[
\pi(m) \, \in \, 
\wedge^2(\fu/(\fl \cap \fu)), \hspace{.2in}
\Pi(\fl) \, \in \, \wedge^2(\fu/\fu_{{\frak l}}).
\]
Thus we only need to show that $\pi(m)$ goes to $\Pi(\fl)$
under the map 
\[
j: \, \fu/(\fl \cap \fu) \lrw \fu/\fu_{{\frak l}}: \, 
x + \fl \cap \fu \Map x + \ful.
\]
But this follows from a general linear algebra fact which we
state as a lemma below.
\qed

\begin{lem}
\label{lem_linear-algebra-fact}
Let $V$ be an $n$-dimensional vector space and let $V^*$ be its 
dual space. On the direct sum vector space $V \oplus V^*$, consider the 
symmetric product $\la\, , \, \ra$ defined by
\[
\la x + \xi, \, y + \eta\ra \, = \, (x, \, \eta) \, + \, (y, \, \xi), 
\hspace{.2in} x, y \in V, \, \xi, \eta \in V^*.
\]

1) Let $V_0$ be any subspace of $V$. For $\lambda \in \wedge^2 (V/V_0)$,
define
\[
W_\lambda \, = \, \{x + \xi: \, x \in V, \, \xi \in V^*, \,
\xi|_{V_0} = 0, \, \xi \backl \lambda = x+V_0.\}
\]
Then $\lambda \mapsto W_{\lambda}$ is  a one-to-one
correspondence between elements  in  $\wedge^2 (V/V_0)$ and maximal
isotropic subspaces $W$ of $V \oplus V^*$ such that $W \cap V = V_0$.

2)  Let $V_1$ be another subspace of $V$ such that $V_0 \subset V_1$.
Let
\[
j: \, V/V_0 \lrw V/V_1: \, v+V_0 \Map v+V_1
\]
be the natural projection. Let $\lambda_0 \in \wedge^2(V/V_0)$
and $\lambda_1 \in \wedge^2(V/V_1)$. Then $j(\lambda_0) = \lambda_1$
if and only if
\begin{equation}
\label{eq_W-lambda-W-jlambda}
W_{\lambda_1} \, = \, V_1 \, + \, (V \oplus V_{1}^{\perp}) \cap
W_{\lambda_0},
\end{equation}
 where $V_{1}^{\perp} = \{\xi \in V^*: \, \xi|_{V_1} = 0
\}.$
\end{lem}

\noindent
{\bf Proof.} 1) Given $\lambda \in \wedge^2(V/V_0)$, it is
easy to see that $W_\lambda$ is maximal isotropic with respect to
$\la \, , \, \ra$ and that $W_\lambda \cap V = V_0$. 
Conversely, if $W$ is a maximal isotropic subspace
of $V \oplus V^*$ such that $W \cap V = V_0$, then
\[
\{\xi \in V^*: \, x + \xi \in W  \, {\rm for} \, {\rm  some} \, x \in V\}
\, = \, V_{0}^{\perp} \, = \,
\{\xi \in V^*: \, \xi|_{V_0} = 0\}.
\]
Define
\[
f: \, (V/V_0)^* \lrw V/V_0: \, \xi \Map x + V_0
\]
where $\xi \in (V/V_0)^* \cong V_{0}^{\perp}$ and
$x \in V$ is such that $x + \xi \in W$. Then $f$ is well defined
and is skew-symmetric. Thus there exists $\lambda \in \wedge^2(V/V_0)$
such that $f(\xi) = \xi \backl \lambda$ for all $\xi \in (V/V_0)^*$.
It is then easy to check that $W = W_{\lambda}$.

2) One way to prove this fact is to take a basis for
$V_0$, extend it first to a basis for $V_1$ and then extend
it further to a basis of $V$. One can then 
write down all the spaces in (\ref{eq_W-lambda-W-jlambda})
using these basis vectors and compare them. We omit the details.
\qed  

As a special case of Theorem \ref{thm_M-to-O}, we have

\begin{cor}
\label{cor_standard-to-canonical}
For any $\fl \in \Lagrfd$ with Property C and any $(U, \piU)$-homogeneous
space $(U/U_1, \pi)$ determined by $\fl$ (see 
Definitions \ref{dfn_C} and \ref{dfn_standard-C}),
the map
\begin{equation}
\label{eq_P}
P: \, (U/U_1, \, \pi) \lrw (U \cdot \fl, \, \Pi): \, 
uU_1 \Map {\rm Ad}_u \fl
\end{equation}
is Poisson. 
\end{cor}

\subsection{Model points}
\label{sec_modelpt}

\begin{dfn}
\label{dfn_model-point}
{\em
We say that a Lagrangian subalgebra $\fl$ is a {\it model point}
(in $\Lagrfd$) if $\fl \cap \fu = \fu_{\frak l}$, where
$\ful$ is the normalizer subalgebra of $\fl$ in $\fu$.
}
\end{dfn}

It is easy to see that the 
set of model points in $\Lagrfd$ is invariant under the $U$-action.

\bigskip
Every model point has Property C, for if $\fl 
\in \Lagrfd$ is a model point,
the connected subgroup
$U_{\frak l}^{'}$ of $U$ with Lie algebra $\fl \cap \fu$
is the identity component of the stabilizer subgroup $U_{\frak l}$
of $\fl$ in $U$, so $U_{\frak l}^{'}$ is closed. Consequently,  
$\fl$ determines a $(U, \piU)$-homogeneous
Poisson structure on any $U/U_1$, where $U_1$ is
a closed subgroup of $U_{\frak l}$, the normalizer subgroup
of $\fl$ in $U$, and has the same Lie algebra $\fl \cap \fu = \fu_{\frak l}$
(see Corollary \ref{cor_Uone} and Definition
\ref{dfn_standard-C}). 
In this
case, the map $\Dr$ in (\ref{eq_P})
is a local diffeomorphism (in addition to being a Poisson map),
and is thus a covering map.
Therefore, the orbit $U \cdot \fl$, together with the
Poisson structure $\Pi$, is a model (up to local diffeomorphism)
of any $(U, \piU)$-homogeneous Poisson space
$(U/U_1, \pi)$ determined by $\fl$. This is the reason we
call $\fl$ a {\it model point} in $\Lagrfd$.

\bigskip
Observe also that $\fl$ is a model point if and only if $T(\fl)=\fl$.

\begin{exam}
\label{exam_model}
{\em
Consider the Lie bialgebra $(\fk, \, \fa + \fn)$ 
in Example \ref{exam_lagrkan}. The 
Lagrangian subalgebra $ \fl = \fa + \fn$ is not a model point
because $\fl \cap \fk = 0$ while the normalizer subalgebra of
$\fl$ in $\fk$ is $\ft = i \fa$. However, $
T(\fl) = \ft + \fn$ is a model point,
as is any real form of $\fg$. In this case, we will show
that every point  in a
certain irreducible component $\Lagr_0$ of $\Lagrfd$ is
a model point. 
}
\end{exam}

When $\fl$ is a 
model point and when its normalizer subgroup $U_{\frak l}$ in $U$
is not connected, the $(U, \piU)$-homogeneous Poisson
spaces $(U/U_1, \pi)$ determined by
$\fl$ might have non-trivial symmetries, as is shown in the following
proposition.

\begin{prop}
\label{prop_cover}
Let $\fl$ be a model point and let $(U/U_1, \pi)$
be any $(U, \piU)$-homogeneous Poisson space 
determined by $\fl$. Then all 
covering transformations for the covering map 
\begin{equation}
\label{eq_covering}
P: \, (U/U_1, \pi) \lrw (U/U_{\frak l}, \Pi): \, uU_1 \Map uU_{\frak l}
\end{equation}
are Poisson isometries for $(U/U_1, \pi)$.
\end{prop}

\noindent
{\bf Proof.} Let $f: U/U_1 \rightarrow U/U_1$
be a covering transformation, so $P \circ f = f$.
We know that $f$ is smooth because it must be of the form
\[
f(uU_1) \, = \, uu_0 U_1
\]
for some $u_0$ in the normalizer subgroup of $U_1$ in $U_{\frak l}$.
Let $x \in U/U_1$
be arbitrary. We need to show that 
$f_* \pi(x) = \pi(f(x))$.
Since $P$ is a local diffeomorphism, it is enough to show that 
$f_* \pi(x)$ and $\pi(f(x))$ have the same image under $P$.
Now since $P$ is a Poisson map and since $P \circ f = f$, we have
\beqa
P_* f_* \pi(x) & = & (P \circ f)_* \pi(x) \, = \, P_* \pi(x) \, = \, 
\Pi(P(x))\\
P_* \pi(f(x)) & = & \Pi (P(f(x))) \, = \, \Pi(P(x)).
\eeqa
Thus $P_* f_* \pi(x) = P_* \pi(f(x))$, and $f$ is a Poisson map.
\qed

In particular, in the case when $U_1 = U_{\frak l}^{'}$ is the identity
connected component of $U_{\frak l}$, 
the group $U_{\frak l} / U_{\frak l}^{'}$
acts on $U/U_{\frak l}^{'}$ as symmetries for 
$(U, \piU)$-homogeneous Poisson structure on $U/U_{\frak l}^{'}$
determined by $\fl$.

\section{Lagrangian subalgebras of $\fg$}
\label{sec_oshima}
\bigskip

In the remainder of the paper, we will concentrate on the Lie bialgebra 
$(\fk, \fa + \fn)$ as described in Example \ref{exam_lagrkan}.
We  first fix more notation.

Throughout the rest of the paper, $\fk$ will be
a compact semi-simple Lie algebra and $\fg = \fk_{\Bbb C}$
its complexification. The Killing form of $\fg$ will be denoted by
$\ll \, , \, \gg$. Let $K$ be a connected  Lie group
with Lie algebra $\fk$ and let $T\subset K$
be a maximal subgroup with Lie algebra $\ft$. Let
$\fh = \ft_{\Bbb C} \subset \fg$ be the complexification
of $\ft$. 
Let $\Sigma$ be the set of
roots of $\fg$ with respect to $\fh$
with the root decomposition
\[
\fg \, = \, \fh \, + \, \sum_{\alpha \in \Sigma} \fg_{\alpha}.
\]
Let $\Sigma_{+}$ be a choice
of positive roots, and let $S(\Sigma_+)$ be the set of simple roots in
$\Sigma_{+}$. We will also say $\alpha > 0$ for $\alpha \in \Sigma_{+}$.
Set $\fa = i \ft$ and let $\fn$ be the complex subspace
spanned by all the positive root vectors. Then  we can
identify $\fk^*$ with $\fa + \fn$
(here
 $\fn$ is regarded as a real Lie subalgebra of $\fg$) through the 
pairing defined by twice the imaginary part of
the Killing form $\ll \, , \, \gg$. This way, $(\fk, \fa + \fn)$
becomes a Lie bialgebra whose double is $\fg = \fk + \fa + \fn$
(Iwasawa Decomposition of $\fg$).  Let $\pi_K$ be the Poisson structure
on $K$ such that $(K, \pi_K)$ is a Poisson Lie group with tangent
Lie bialgebra $(\fk, \fa + \fn)$. We can describe $\pi_K$ 
explicitly as follows:
Let $\theta$ be the complex conjugation of $\fg$ defined by $\fk$.
Let $\ll \, , \, \gg_{\theta}$ be
the Hermitian positive definite inner product on $\fg$ given by
\[
\ll x, y\gg_{\theta} \, = \, - \ll x, \, \theta y\gg, \hspace{.2in}
x, y \in \fg.
\]
For each $\alpha \in \Sigma_{+}$, choose $\ea \in \fg_{\alpha}$ such that
\[
\ll \ea, \, \ea \gg_{\theta} \, = \, 1.
\]
Let $\eb = -\theta(\ea) \in \fg_{-\alpha}$ so that $\ll \ea, \, \eb \gg
= 1$. Set
\[
\Xa \, = \, \ea - \eb \, = \, \ea + \theta(\ea), \hspace{.3in}
\Ya \, = \, i(\ea + \eb) \, = \, i\ea + \theta(i\ea).
\]
Then
\[
\fk \, = \, \ft \, + \, {\rm span}_{\Bbb R}\{\Xa, \Ya: \alpha \in \Sigma_{+}\}.
\]
The Poisson bivector field on $K$ is 
given by
\[
\pi_K(k) \, = \, r_k \Lambda \, - \, l_k \Lambda, \hspace{.2in}
k \in K,
\]
where
\[
\Lambda \, = \, {1 \over 4} \sum_{\alpha \in \Sigma_{+}} \Xa
\wedge \Ya \, \in \, \fk \wedge \fk.
\]

Recall that 
a real
subalgebra $\fl$ of $\fg$ is {\it Lagrangian} if
${\rm Im} \ll x, \,  y\gg = 0$ for all $x, y \in \fl$ and if
$\dim_{\Bbb R} \fl = \dim_{\Bbb C} \fg$.
These Lagrangian subalgebras correspond to $(K,\pi_K)$
Poisson-homogeneous spaces by Drinfeld's theorem. The 
set of all Lagrangian subalgebras of $\fg$ will be denoted by
$\Lagr$. 
It is an algebraic subset of the Grassmannian ${\rm Gr}(n, \fg)$ of
$n$-dimensional subspaces of $\fg$ (regarded as a $2n$-dimensional
real vector space). 

In this section,  we will decompose ${\cal L}$ into a finite union
of manifolds.

\subsection{Karolinsky's classification}
\label{sec_kara-classification}

\bigskip
E. Karolinsky \cite{karo:homog-compact} 
has determined all Lagrangian subalgebras
$\fl$ of $\fg$. 
To describe his result, we need some notation. Let
$S\subset S(\Sigma_+)$ be a subset of the set of simple roots, and
let $[S]$ be the set of roots in the linear span of $S.$
Consider 
\[
\fmS = \fh \oplus ({\bigoplus}_{\alpha \in [{\scriptscriptstyle S}]}^{} 
{\fg}_\alpha),
\hspace{.2in}
\fnS = {\bigoplus}_{\alpha\in {\Sigma}_+ -[{\scriptscriptstyle S}]}^{}
 {\fg}_\alpha
\]
and 
\[
\fpS = \fmS + \fnS,
\]
so that $\fpS$ is a parabolic subalgebra of type $S,$ $\fnS$
is its nilradical, and $\fmS$ is a Levi factor. Let
$\fmSss = [\fmS, \fmS]$ be the (semi-simple) derived algebra
of $\fmS$. The center 
of $\fmS$ is
\begin{equation}
\label{eq_zs}
\fzS = \{ H\in \fh: \,   \alpha_i(H)=0, \, \forall \alpha_i\in S \},
\end{equation}
which is also the orthogonal complement of $\fmSss$ in
$\fmS$ with respect to the Killing form of $\fg$ restricted to
$\fmS$.
Thus the restriction of the  Killing form to $\fzS$ is nondegenerate,
and we may consider Lagrangian subspaces of $\fzS$ (regarded
as a real vector space) with respect
to the restriction to $\fzS$ of the imaginary part of the Killing form.

\bigskip
Now for any subset $S$ of the set of simple roots, a Lagrangian subspace
$V$ of $\fzS,$ and a real form 
${\fm}_{\scriptscriptstyle S, 1}^\tau$ of $\fmSss,$ set 
\[
\fl(S,V,\tau)={\fm}_{\scriptscriptstyle S, 1}^\tau \oplus V \oplus \fnS.
\]
It is easy to see that it is a Lagrangian subalgebra of $\fg$.

\begin{dfn}
\label{dfn_stdlag}
{\em
We will call $\fl(S,V,\tau)$ the  {\it standard Lagrangian subalgebra
associated to $(S,V,\tau)$}.
}
\end{dfn}

\begin{thm}
\cite{karo:homog-compact}
\label{thm_karolinsky}
{\it 
Every Lagrangian subalgebra of $\fg$ is of the form ${\rm Ad}_k
(\fl(S,V,\tau))$
for some $k\in K.$
}
\end{thm}

Note that the nilradical of ${\rm Ad}_k(\fl(S,V,\tau))$ is ${\rm Ad}_k(\fnS)$.
Denote by $\PS$ the connected subgroup of $G$ with Lie algebra $\fpS$.

\begin{prop}
\label{prop_uniqueness}
{\it
Let 
\[
\fl = {\rm Ad}_k (\fl(S,V,\tau)) = {\rm Ad}_{k_1} (\fl(S_1,V_1,\tau_1))
\]
 be a Lagrangian subalgebra. Then $S=S_1$, $V=V_1$, $k^{-1}k_1\in \PS$,
and $\tau$ is conjugate to $\tau_1$ in $K\cap \PS$.
}
\end{prop}

\noindent
{\bf Proof.}
We have ${\rm Ad}_{k^{-1}k_1} (\fl(S_1,V_1,\tau_1))=\fl(S,V,\tau)$.
Using the fact that conjugate algebras have conjugate nilradicals, it
follows easily that ${\rm Ad}_{k^{-1}k_1} \fnSone=\fnS$.
 From the
definition of $\fnS$, it follows that $S=S_1$. Moreover, since
$\fpS$ is the perpendicular complement of $\fnS$, it
 follows that ${\rm Ad}_{k^{-1}k_1}$ normalizes $\fpS$. Since a
parabolic subgroup is the normalizer of its nilradical, $k^{-1}k_1 \in
\PS$. The remaining claims follow from the facts that $\fnS$ is an
ideal and $\fzS$ is central in $\fmS$.
\qed

%

In the following, we study separately the pieces that come into the Karolinsky
classification. 

\subsection{Lagrangian subspaces of $\fzS$}
\label{sec_lzs}
For a subset $S$ of the set of simple roots, let $\fzS$
be given as in (\ref{eq_zs}).
Since the Killing form is nondegenerate on $\fzS$, 
its imaginary part $B$ is a nondegenerate
symmetric bilinear form of index $(z,z)$ on $\fzS$, now regarded
as  a $2z$-dimensional real vector space. Denote by 
$\LagrfzS$ the variety of Lagrangian subspaces of $\fzS$ with
respect to $B.$

\begin{prop}
\label{prop_indlagr}
The variety
\[
\LagrfzS =\cup_{\epsilon = \pm 1} \LagrfzSeps
\]
 is a
smooth manifold of dimension ${{z(z-1)}\over 2}$
with two connected components $\LagrfzSeps$, $\eps = \pm 1 $. We call
$\LagrfzSplus$  the component containing $\fzS \cap \ft$ and
call $\LagrfzSminus$ the other one. Each component is Zariski closed.
\end{prop}

\noindent
{\bf Proof.} 
The first assertion follows from the identification of $\LagrfzS$ with
$O(n)\times O(n)/O(n)$ given in \cite{port}, Theorem 14.10.
The algebraicity of each of the components can be derived from
the discussion of charts in \cite{port} following Theorem 14.10,
or by noting the corresponding fact for the space $\LagrfzSC$ of 
complex linear Lagrangian subspaces of the complexification $\fzSC$
with respect to the nondegenerate Killing form
(see \cite{ACGH}, Exercise B, pp. 102-103),
and verifying the easy fact that $\LagrfzS$ is the set of real
points of $\LagrfzSC$. 
\qed

We remark that two Lagrangian subspaces $V$ and $V^\prime$
lie in the same component if and only if $\dim (V\cap V^\prime) = \dim (V)
{\rm mod } \ 2$. This is proved in the complex case in \cite{ACGH},
and the real case can be deduced from the complex case. It follows
that $\ft\cap\fzS$ and $\fa\cap\fzS$ lie in the same component if and
only if $\dim (\fzS)$ is even.

\subsection{Real forms of  $\fg$}
\label{sec_real}

A real form of $\fg$ is clearly a Lagrangian subalgebra of $\fg$.
Denote by $\cR$ the set of all real forms of $\fg$. We will 
recall some facts about
$\cR$ in this section (see \cite{onitvin} or \cite{abv} for more
details.) 

Let
$\As$ be the group of complex linear
automorphisms of $\fg.$ Its identity component is the adjoint group
$G= \Ints$ of interior automorphisms of $\fg$. Let
$\Outs$ be the automorphism group of the Dynkin diagram of $\fg$.
It is well-known that there is a split short exact sequence
\[
0 \, \longrightarrow \Ints \, \longrightarrow \, \As \, 
\stackrel{\phi}{\longrightarrow} \, \Outs \, \longrightarrow \, 0.
\]

Let $\theta$ be the Cartan involution of $\fg$ defined by
the compact real form $\fk$. We will identify a real form $\fg_0$ 
of $\fg$ with the complex conjugation $\tau$
on $\fg$ such that $\fg_0=\fg^\tau$. Define a map
\[
\psi: \, \cR \, \longrightarrow \, \Outs
\]
as follows:
\[
\psi(\tau) \, = \, \phi(\tau \theta) \, = \, \phi(\theta \tau).
\]
To see that $\phi(\tau \theta) \, = \, \phi(\theta \tau)$, choose
$g \in \Ints$ be such that $\tau_1 = g \tau g^{-1}$ commutes
with $\theta$ (see \cite{hel}, Theorem III.7.1, and the following
remark). Then
 we get
\[
\phi(\tau \theta) = \phi(g^{-1} \tau_1 g \theta) = \phi(\tau_1 g \theta) = 
\phi(\tau_1 \theta \theta^{-1} g \theta) = \phi(\tau_1 \theta)
\]
and similarly, $\phi(\theta \tau)=\phi(\theta \tau_1)$.
Since $\tau_1$ commutes with $\theta$, we have 
$\phi(\tau \theta) \, = \, \phi(\theta \tau)$.
In particular, we see that $\psi(\tau)$ is an involution.

Conversely, let $d$ be an involutory automorphism of the Dynkin diagram
$D(\fg).$ Then $d$ extends to a complex linear involution $\dgamma$ of
$\fg$ as follows: we can choose $\dgamma \in \As$ preserving $\fh$
and permuting the fixed simple root vectors $E_\alpha, \alpha \in S(\Sigma_+)$
(see for example the proof of Proposition 2.7 in \cite{abv}).
Then $\dgamma(E_{\alpha})=E_{d\alpha}$ and $\dgamma(E_{-\alpha})=
E_{-d\alpha}$. If $H_\alpha=[E_\alpha,E_{-\alpha}],$ it follows
that $\dgamma(H_\alpha)=H_{d\alpha}$, and also $\dgamma$ commutes with the
Cartan involution on generators, and therefore on all of $\fg$.

Set
\[
{\cal L}(\fg, d) \, = \, \psi^{-1}(d).
\]
Then
\[
{\cR} \, = \, \cup_d \Lagr(\fg, \daut)
\]
is a finite disjoint union, where $d$ runs over the set of all
involutory diagram automorphisms of $\fg$. 

Let 
$\dtheta = \dgamma \theta = \theta \dgamma$.
Then $\tau_d \in \Lagr(\fg, d)$. To describe all the elements in 
$\Lagr(\fg, d)$, consider
\[
\Gmin 
      =\{ g\in \Ints: \, (g \dtheta)^2 = 1 \}
=\{ g\in \Ints: \, \dtheta(g)=g^{-1} \}.
\]
If $g\in \Gmin$, then
$g \dtheta$ is a real form of $\fg$ and
$\psi(g \tau_d) = d$, so $g \tau_d \in \Lagr(\fg, d)$.
Conversely, if $\tau \in {\cal L}(\fg, d)$, then
$\phi(\tau \theta) = \phi(\dgamma)$, so $\tau = g \dtheta$ for
some $g \in \Ints = \ker(\phi)$. But $\tau^2 = 1$, so
$g \in \Gmin$. Hence every real form $\tau$ in ${\cal L}(\fg, d)$ is 
of the form $\tau = g \dtheta$ for some $g \in \Gmin$.

\begin{lem}
\label{lem_selfnorm}
Every real form of $\fg$ is its own normalizer in $\fg$.
\end{lem}

\bigskip
\noindent {\bf Proof.} The proof follows easily by considering the $\pm 1$ 
eigenspace decomposition 
$\fg={\fg}^{\tau}\oplus {\fg}^{-\tau}$ of $\tau$. 
\qed

\begin{lem}
\label{lem_lagrsemisimple}
$\Lagr(\fg,\daut)$ is a smooth submanifold of
$\Grg$ of  dimension $\dim_{\Bbb C} \fg$.
\end{lem}

\noindent
{\bf Proof.} 
Note that $\Ints$ acts on $\Lagr(\fg,\daut)$ by the action
$g\cdot \tau =g \tau g^{-1}.$ The orbits of the larger group
$\As$ on  the set of all real forms are the equivalence classes
of real forms, and there are only finitely many of them (see \cite{onitvin}).
Since $\Ints$ is the identity connected component of $\As$
and $\As$ has only finitely many components, it follows that
$\Ints$ has only finitely many orbits on the set of all real forms.
Since $\Lagr(\fg,\daut)$ is a subset of the set of all real forms,
it follows that $\Lagr(\fg,\daut)$ is a finite union of $\Ints$ orbits.
Now the action of $\Ints$ on ${\rm Gr}(n, \fg)$ by $(g, \fl) \mapsto g(\fl)$ 
is smooth and ${\cal L}(\fg, d) \subset {\rm Gr}(n, \fg)$ is a disjoint
union of finitely many $\Ints$-orbits, it follows that each 
$\Ints$-orbit in ${\cal L}(\fg,d)$ is a smooth submanifold
of ${\rm Gr}(n, \fg)$. Moreover, by Lemma 
\ref{lem_selfnorm}, 
all orbits have the same dimension.
Thus, $\Lagr(\fg,\daut)$ is a smooth submanifold of ${\rm Gr}(n, \fg)$
of dimension $\dim_{\Bbb C} \fg$.
\qed

We will show later that the closure of $\Lagr(\fg,\daut)$ in
$\Lagr$ is a smooth, compact and  connected submanifold of
${\rm Gr } (n,\fg)$.

\subsection{Model points}
\label{sec_modpts}

\begin{lem}
\label{lem_lagrnormalizer}
The normalizer of the Lagrangian subalgebra $Ad_k(\fl(S,V,\tau))$
in $\fg$ is 
\[
Ad_k(\fr(S,\tau)):= 
Ad_k({\fm}_{\scriptscriptstyle S, 1}^\tau \oplus \fzS \oplus \fnS).
\]
\end{lem}

\bigskip
\noindent {\bf Proof.} It suffices to prove the statement when 
$k=e$, the identity element of $K$. 
It is clear that $\fr(S,\tau)$ normalizes $\fl(S,V,\tau)$.
Conversely, 
if $X\in \fg$ normalizes $\fl(S,V,\tau),$ it normalizes
its nilradical $\fnS,$ so 
it normalizes the perpendicular 
$\fpS$ of $\fnS.$ Since $\fpS$ is
parabolic, it equals  its own normalizer, so $X\in \fpS.$
Write $X=X_1+X_2,$ with $X_1\in \fmS$ and  $X_2\in \fnS.$
Then $X_1$ normalizes ${\fm}_{\scriptscriptstyle S, 1}^\tau$.
It follows from Lemma \ref{lem_selfnorm} that
 $X_1\in {\fm}_{\scriptscriptstyle S, 1}^\tau  +\fzS.$
\qed

\begin{prop}
\label{prop_modelpt}
The Lagrangian subalgebra $Ad_k(\fl(S,V,\tau))$ is a model point
if and only if $V=\fzS \cap \ft$.
\end{prop}

\bigskip
\noindent {\bf Proof.} Since the set of model points is $K$-invariant,
 it suffices to prove the proposition when
$k=e$. Let $N_{\frak k}(\fl(S,V,\tau)$ be the normalizer
of $\fl(S, V, \tau)$ in $\fk$. By the previous lemma, the
quotient
\[
N_{\frak k}(\fl(S,V,\tau)/\fk\cap \fl(S,V,\tau) = (\fzS \cap \ft)/(V\cap \ft),
\]
since $\fzS\cap \fk = \fzS\cap \ft$. 
The proposition now  follows from the definition of model points.
\qed

\begin{rem}
\label{rem_Tsq}
{\em
In fact, 
essentially the same argument shows that if $\fl=Ad_k(\fl(S,V,\tau))$,
then $T(\fl)=Ad_k(\fl(S,\fzS\cap \ft,\tau))$ (see 
Theorem \ref{thm_Pi-Lagr-at-fl}
for the definition of $T(\fl)$). It follows that $T(T(\fl))=
T(\fl)$ for $\fl \in \Lagr$. For a general Lie bialgebra, $T\circ T \not=
T$. Indeed, for a Lie algebra $\fu$, we can form a Lie bialgebra
$(\fu, \fu^*)$, where  ${\fu}^*$  has the abelian Lie algebra structure.
Its double is the semi-direct product Lie algebra structure on
$\fu + \fu^*$ defined by the co-adjoint action of $\fu$ on 
$\fu^*$. 
Consider the case when $\fu$ is the three dimensional Heisenberg
algebra with basis $\{X,Y,Z \}$ with $Z$ central and $[X,Y]=Z$, and
let $f_X,f_Y,f_Z$ be the dual basis. Let $\fl$ be the Lagrangian
subalgebra spanned by $X,f_Y$ and $f_Z$. Then $T(\fl)$ is spanned by
$X, Z$ and $f_Y$ while $T(T(\fl))=\fu$.
}
\end{rem}

\begin{cor}
\label{cor_Gmodelpts}
$G$ preserves the set of model points.
\end{cor}

\bigskip
\noindent{\bf Proof.} It suffices to consider model points 
$\fl(S,\fzS \cap\ft,\tau)$. Let $\PS, \MS$ and  $\NS$ be the connected
Lie groups with Lie
algebra $\fpS$, $\fmS$ and $\fnS$ respectively.
 Since $K$ acts transitively on $G/\PS$ and preserves the set
of model points, it
suffices to prove that $Ad_p(\fl(S,\fzS\cap\ft,\tau))$ is a model
point for $p\in \PS$. Using the Levi decomposition $\PS=\MS \NS$
we write $p=mn$. Since 
$Ad_n\fl(S,\fzS\cap\ft,\tau)=\fl(S,\fzS\cap\ft,\tau)$, it suffices
to prove that $Ad_m$ preserves model points in $\fpS$, which follows because
$M$ acts trivially on $\fzS$.
\qed

\begin{rem}
\label{rem_Gnotmodelpts}
{\em
In general, the adjoint group of the double Lie algebra does
not preserve the set of model points. Indeed, let
$\fg$ be a semisimple Lie algebra with triangular decomposition
$\fg=\fn+\fh+\fnm$, Borel subalgebra $\fbp=\fh+\fn$ and opposite
Borel $\fbm=\fh+\fnm$. Then the Lie algebra $\fd=\fg\oplus \fh$ is the
double of the pair $(\fbp,\fbm)$ with embeddings $i_{\pm}:\fbpm\to \fd$
given by $i_{\pm}(H+x)=(H+x,\pm H)$ with $H\in \fh$, $x\in \fn$ or $\fnm$.
Let $n\in N_G(\ft)$ be a representative for the long element of
the Weyl group. Then although $\fbp$ is clearly a model point,
$Ad_n (\fbp)$ is not a model point.
}
\end{rem}

\subsection{Lagrangian data}
\label{sec_lag-data}

\begin{dfn}
\label{dfn_lagdata}
{\em 
A triple $(S,\eps,\daut)$ is called Lagrangian datum if $S\subset S(\Sigma_+)$
is a subset of the set of simple roots, $\eps = \pm 1$, and $\daut$
 is a diagram
automorphism for the Dynkin diagram $D(\fmSss)$ of $\fmSss$. If 
$\fl ={\rm Ad}_k(\fl(S,V,\tau))$, $k\in K,$ is
a Lagrangian subalgebra, then $\fl$
 has associated Lagrangian
data $\Phi(\fl)=(S,\eps,d)$,
where
$\eps=1$ if $V$ lies in the same 
connected component of $\LagrfzS$ as $\fzS
\cap \ft$ and is $-1$ otherwise,
and $\daut$ is the diagram automorphism  of $\fmSss$
defined by $\tau.$ It follows from Proposition \ref{prop_uniqueness}
that the triple $(S,\eps,d)$ is determined by $\fl$.
}
\end{dfn}

Given Lagrangian datum $(S,\eps,\daut),$ we let
\[
\Lagr(S,\eps,\daut) = \{ \fl: \,  \Phi(\fl)=(S,\eps,\daut) \}.
\]
Then
\[
{\cal L} \, = \, \cup_{({\scriptscriptstyle S}, \epsilon, d)} 
\Lagr(S,\eps,\daut).
\]
Note that this is a finite disjoint union.

\begin{prop}
\label{prop_lsedfiber}
For each Lagrangian datum $(S, \epe, d)$, 
$\Lagr(S,\eps,\daut)$ is a smooth submanifold of the 
Grassmannian ${\rm Gr}(n, \fg)$ of dimension
$\dim(\fk)+{{z(z-3)}\over 2}$, and it fibers over $G/\PS$ with
the fiber being the product of $\LagrfzSeps$ and $\Lagr(\fmSss,\daut)$.
\end{prop}

\noindent
{\bf Proof.} 
Consider the subset 
\[
\Lagr_{\fpS} (S,\eps,\daut) = \{ \fl(S,V,\tau) \}
\]
of all standard Lagrangian subalgebras (see Definition \ref{dfn_stdlag})
attached 
to the Lagrangian datum $(S,\eps,d)$. 
It can be identified with
$\Lagr(\fmSss,\daut)\times  \LagrfzSeps$ as a submanifold of the
Grassmannian ${\rm Gr} (n,\fg)$.
 Indeed, $\Lagr(\fmSss,\daut)$ is a submanifold of the
Grassmannian of ${\rm Gr} (m,\fmSss)$
 where $m=\dim (\fmSss)$, $\LagrfzSeps$ is
a submanifold of ${\rm Gr}(z,\fzS)$, and the direct sum map 
${\rm Gr} (m,\fmSss)\times {\rm Gr} (z,\fzS) \to {\rm Gr }(n,\fg)$,
 $(U,V)\mapsto 
U\oplus V\oplus \fnS$ is a closed embedding.

We consider the multiplication map
\[
m:K\times _{K\cap P_S} \Lagr_{\fpS} (S,\eps,\daut) \to
\Lagr (S,\eps,\daut), \ \ m(k,\fl ) = Ad_k(\fl)
\]
The fiber product is a smooth manifold since it is a fiber bundle
over $K/K\cap \PS\cong G/\PS$ with smooth fiber $\Lagr_{\fpS} (S,\eps,\daut)$.
The map $m$
 is onto by the Karolinsky classification Theorem \ref{thm_karolinsky},
and it is clearly smooth and proper. We will show that it is an immersion,
and it will follow that ${\Lagr (S,\eps,\daut)}$ is a 
smooth submanifold of ${\rm Gr } (n,\fg)$. 

The fact that $m$ is injective follows from Proposition 
\ref{prop_uniqueness}. In order to show that the tangent map
$m_*$ is injective, it suffices to show $m_*$ is injective
at points of the form $(e,\fl(S,V,\tau))$ by $K$-equivariance.
Recall that the tangent space at a plane $U$ to the Grassmannian
${\rm Gr}(n,V)$ of $n$-planes in a space $V$
 can be identified with
$\Hom (U,V/U)$.
 Using this identification, the tangent space
to the fiber product 
$K\times_{K\cap P_S} {\rm Gr}(n,\fpS)$ at $\fl(S,V,\tau)$ 
is the quotient of $\fk \oplus \Hom (\fl(S,V,\tau),
\fpS/\fl(S,V,\tau))$ by the relation $(X-Y,\xi(Y)+Z)\sim (X,Z)$, where
$X\in \fk$, $Y\in \fk \cap \fpS,$ $\xi(Y)$ is the induced vector
field at $\fl(S,V,\tau)$, and $Z\in \Hom (\fl(S,V,\tau),\fpS/\fl(S,V,\tau))$.
Observe that for $Z$ to be tangent to the fiber $\Lagr_{\fpS}(S,\eps,d)$,
we must have $Z:\fnS\to 0$. 
When we identify the tangent space to ${\rm Gr}(n,\fg)$ 
at $\fl(S,V,\tau)$ with
$\Hom(\fl(S,V,\tau),\fg/\fl(S,V,\tau))$,  the tangent map
is $m_*(X,Z)=\xi(X)+Z$, where $\xi(X)$ is the induced vector
field. Now the claim that $m_*$ is injective follows since for
any $X\not\in \fk\cap\fpS,$ $\xi(X)\cdot \fnS\not\subset \fl(S,V,\tau)$.
To verify this last assertion, let $X\in \fk\backslash\fk\cap \fpS$,
and choose a maximal root $\alpha \not\in [S]$ such that the projection
$p_{-\alpha}(X)$ of $X$ to the ${\fg}_{-\alpha}$ root space is nonzero.
Then $[X,{\fg}_\alpha]=[p_{-\alpha}(X),{\fg}_\alpha] + Y$ where
$\ll Y, Y \gg = \ll Y, [p_{-\alpha}(X),{\fg}_\alpha] \gg = 0$.
Since $[p_{-\alpha}(X),{\fg}_\alpha]= [{\fg}_\alpha,{\fg}_{-\alpha}]$,
which is a $2$-dimensional real vector space on which the imaginary
part of the Killing form is not isotropic, 
it follows that $[X,{\fg}_\alpha]$ is not isotropic. Thus, $[X,{\fg}_\alpha]$
is not contained in any Lagrangian subalgebra.

The dimension statement follows from Proposition \ref{prop_indlagr}
and Lemma \ref{lem_lagrsemisimple}. 
\qed

\begin{rem}
\label{rem_Glsed}
{\em
Note that $G$ preserves $\Lagr(S,\eps,d)$. The proof is similar
to that of Corollary \ref{cor_Gmodelpts}.
}
\end{rem}


\begin{exam}
\label{exam_lsed}
{\em
When $S$ is the set of all simple roots, we have $\fmS = \fg$ and $\epsilon$ 
can only be $1$, so ${{\cal L}}(S, \epsilon, d) = {{\cal L}}(\fg, d)$.
}
\end{exam}
\qed

\begin{exam}
\label{exam_morelsed}
{\em  For $\fg=\fsl (2, \C)$, there are three $\Lagr(S,\eps,\daut)$'s.
First, 
$\Lagr(S(\Sigma_+),1,{\rm id})$ is a disjoint union of the two symmetric
spaces $SO(3,\C)/SO(3,\R)$ and $SO(3,\C)/SO(2,1)$, where the first
piece consists of compact real forms and the second piece consists
of real forms isomorphic to $\fs\fo (2,1)$. $\Lagr(\emptyset,1,
{\rm id})$
is the $SL(2,\C)$ orbit of $\ft+\fn$ and is isomorphic to
$\C P^1$. $\Lagr(\emptyset,-1,{\rm id})$ is also isomorphic to $\C P^1$,
and is the $SL(2,\C)$ orbit through $\fa+\fn$.
As we will show in Section \ref{sec_lsedclos}, $\Lagr(\emptyset,1,{\rm id})
\subset \overline{\Lagr(S(\Sigma_+),1,{\rm id})}$. 
This last closure can be identified
with $\R P^3$, the projectivization of $2\times 2$ Hermitian matrices.

In case $\fg=\fsl (3)$, there are eight $\Lagr(S,\eps,\daut)$'s.
$\Lagr(S(\Sigma_+),1,{\rm id})$ is a union of components consisting of the
real forms isomorphic to $\fs\fu (p,3-p)$. It is a union of
symmetric spaces. Let $\sigma$ be the nontrivial involution
of the Dynkin diagram of $\fsl (3)$. 
Then $\Lagr(S(\Sigma_+),1,\sigma)$ consists of real
forms isomorphic to $\fsl (3,\R)$. There are four pieces of the
form $\Lagr(\alpha_i,\pm 1,{\rm id})$ corresponding to the two choices
of $\alpha_i$ and the two choices of $\pm 1$. Each of these
pieces fibers over $G/P_i$ for a parabolic $P_i$ with  the fiber being a
symmetric space for $SL(2,\C)$. The final two components are
of the form $\Lagr(\emptyset,\pm 1,{\rm id})$. These are bundles over
the full flag variety $G/B$ with the fiber being a component of the variety
of Lagrangian subspaces of $\R ^4$ with respect to a quadratic form
of index $(2,2)$. The only nontrivial inclusions are $\Lagr(\alpha_i,1,
{\rm id})\subset \overline{\Lagr(S(\Sigma_+),1,{\rm id})}$.
}
\end{exam}

Because of the fiber bundle decomposition of $\Lagr(S,\eps,\daut)$ and the
fact that the base and $\LagrfzSeps$ are compact, the
study of the closure $\overline{\Lagr(S,\eps,\daut)}$ can be
reduced to the study of $\overline{\Lagr(\fg, \daut)}$ for $\fg$
semisimple.
In the  following Sections \ref{sec_sig} and \ref{sec_dp},
we show that
$\overline{\Lagr(\fg, \daut)}$  is a smooth connected
submanifold of $\Grg$. We will also determine
its decomposition into $G$-orbits. 
These results will be applied in Section \ref{sec_lsedclos}
to show that $\overline{\Lagr(S,\eps,\daut)}$ is a smooth
submanifold of ${\rm Gr }(n,\fg)$.

\section{Extended signatures and the corresponding Lagrangian subalgebras
of $\fg$}
\label{sec_sig}

In this section, we give examples of Lagrangian subalgebras of $\fg$ that lie
in $\overline{\Lagr(\fg,d)}$. They are 
obtained by
considering 
{\it extended signatures} of roots of $\fg$ as slightly
generalized from \cite{os:compact}. They will be used in Section \ref{sec_dp}
to describe $G$-orbits in $\overline{\Lagr(\fg,d)}$.

\subsection{Extended signatures}
\label{subsec_extsig}

Recall that 
\[
S(\Sigma_+) \, = \, \{\alpha_1, \alpha_2, ..., \alpha_l\}
\]
is the set of simple roots in $\Sigma_{+}$.
Let $d$ be an involutory automorphism of the Dynkin diagram of $\fg$.

\begin{dfn}
\label{dfn_sig}
{\em 
An  {\it extended
$\daut$-signature} of the root system $\Sigma$
is a map $\sigma:\Sigma\to \{ -1,0,1 \}$ satisfying
\begin{equation}
\label{eq_sign1}
\sigma(\alpha)=\prod {\sigma(\alpha_i)}^{m_i}, \ \rm{where } \
\alpha = \sum_{i=1,\dots, l} m_i\alpha_i
\end{equation}
\begin{equation}
\label{eq_sign2}
\sigma(\daut(\alpha_i))=\sigma(\alpha_i).
\end{equation}
We say that $\sigma$ is a {\it $\daut$-signature}
 if $\sigma(\alpha) \neq 0$ for any
$\alpha \in \Sigma$. 
}
\end{dfn}

An extended $d$-signature $\sigma$
 is determined by its value on the simple roots.
If $\sigma$ is an extended $d$-signature,
 let
\[
{\rm supp}(\sigma)= \{ \alpha \in \Sigma: \sigma(\alpha)\not= 0 \}.
\]
Then $S_{\sigma}: = 
S(\Sigma_+) \cap \supp (\sigma)$ is $\daut$-invariant. If we use
$[S_{\sigma}]$ to denote the set of roots that
are  in the linear span of $S_{\sigma}$, then
\[
\supp (\sigma) \, = \, [S_{\sigma}].
\]
Let
\[
S_{\sigma,1} \, = \, \{\alpha_i \in S(\Sigma_+): \, \sigma(\alpha_i) = -1\},
\hspace{.2in}
\rhocheckone \, = \, 
\sum_{\alpha_i \in S_{\sigma,1}} {\check{h}}_i \, \in \, \fa,
\]
where $\{{\check{h}}_i: i = 1, .., l\} \subset \fa$ 
is the set of fundamental coweights
corresponding to the simple roots, namely
$\alpha_i({\check{h_j}}) = \delta_{i,j}$ for $i, j = 1, ..., l$.
Then
\begin{equation}
\label{eq_userho}
\sigma(\alpha) \, = \, \left\{ \begin{array}{ll} 0, & \alpha \notin
[S_{\sigma}] \\
(-1)^{\alpha(\rhocheckone)}, & \alpha \in [S_{\sigma}].
\end{array} \right.
\end{equation}
Conversely, for any $\daut$-invariant subset 
$S$ of $S(\Sigma_+)$ and any $\daut$-invariant subset $S_1$ of $S$,
there is an extended $d$-signature $\sigma$ such that $S = S_{\sigma}$
and $S_1 = S_{\sigma,1}$.

\bigskip
For an extended $d$-signature $\sigma,$ let 
\[
\fm_\sigma = {\fm}_{{\scriptscriptstyle S}_\sigma} = \fh \oplus 
\left(\bigoplus_{\alpha \in [S_{\sigma}]} \fg_{\alpha}\right),
\, \,  \fn_\sigma = {\fn}_{{\scriptscriptstyle S}_\sigma} = 
\bigoplus_{\alpha \in \Sigma_{+} - [S_{\sigma}]} \fg_{\alpha}, 
\, \, 
\fp_{\sigma} = \fp_{{\scriptscriptstyle S}_\sigma} = \fm_\sigma \oplus 
\fn_\sigma
\]
as in the notation in Section \ref{sec_kara-classification}. Also let
$\fz_\sigma = \fz_{{\scriptscriptstyle S}_\sigma}$ 
be the center of $\fm_\sigma$, and let
\[
{\fn_\sigma}_-=\bigoplus_{\alpha \in -{\Sigma}_+,  \sigma(\alpha)=0 } {\fg}_\alpha,
\hspace{.2in}
\fmSigss = [\fm_\sigma,\fm_\sigma].
\]
 Then $\sigma$ determines 
a complex linear involution $a_\sigma$ of $\fm_{\sigma}$ by 
\[
 a_\sigma |_{\frak h} = {\rm id}, \, \, \, 
a_\sigma |_{{\frak g}_\alpha} = 
\sigma(\alpha)\cdot {\rm id},
\]
where $\alpha \in {\rm supp}(\sigma)$.
In other words,
\[
a_{\sigma} \, = \, {\rm Ad}_{\exp(\pi i \rhocheckone)}.
\]
Let $\dtheta=\dgamma \theta$ be the conjugate linear involution of
$\fg$ discussed in Section \ref{sec_real}. 
Then it is routine to check that $\tau_{d,\sigma} := a_\sigma 
\dtheta$
is a conjugate linear involution of $\fm_\sigma$ so the Lie algebra
\[
    {\fk}_{d,\sigma} = \fm_{\sigma}^{\tau_{d,\sigma}}
\]
is a real form of $\fm_\sigma.$ 
Set
\[
\fl_{d,\sigma}={\fk}_{d,\sigma}+\fnsigmaplus.
\]
It is easy to check that $\fls$ is a Lagrangian subalgebra of $\fg$.

Since $S_{\sigma}$ is $d$-invariant, $\fm_\sigma$ is invariant
under $\gamma_d$. Regarded as an complex automorphism of $\fmSigss$, 
$\gamma_d$ defines an automorphism of the Dynkin diagram
of $\fmSigss$ which is just $d|_{S_\sigma}$. Let $\fz_\sigma^{\tau_d}$ be the
fixed point set of $\tau_d$ restricted to $\fz_\sigma$. Set $\epsilon = 1$
if $\fz_\sigma^{\tau_d}$ lies in the same component as $ 
\fz_\sigma \cap \ft$ and $\epsilon = -1$
otherwise. Then, since $a_\sigma$ is an inner automorphism
of $\fmSigss$, we know that
 $\fl_{d, \sigma} \in \Lagr(S_\sigma, \epsilon, d|_{S_\sigma})$.

\begin{exam}
\label{exam_flds}
{\em
When $\sigma(\alpha) = 0$ for all $\alpha$, we have $\fl_{d, \sigma}
= \fh^{\tau_d} + \fn$. On the other hand, $\sigma$ is a $d$-signature if
and only if $\flds$ is a real form of $\fg$. In this case,
$\flds \in \Lagr(\fg, d)$.
}
\end{exam}


We choose $H\in \fh$ such that $\alpha(H) > 0$ if $\alpha$ is
a root of ${\fn}_\sigma$ and $\alpha(H) = 0$ if $\alpha$
is a root of $\fmSigss$. Choose a $\daut$-signature $\sigma^\prime$
such that $\sigma^\prime (\alpha)=\sigma(\alpha)$ if $\sigma(\alpha)
\not= 0$. Then by writing down generators of ${\fl}_{d,\sigma^\prime}$,
one can check that 
\[
\lim_{t\to +\infty} \exp(tH){\fl}_{d,\sigma^\prime}={\fl}_{d,\sigma},
\]
where the limit takes place in the Grassmannian $\Grg$.
It follows that ${\fl}_{d,\sigma}\in \overline{\Lagr(\fg,\daut)}$.

\subsection{Extended signatures and real forms}
\label{subsec_realforms}

\bigskip
To relate real forms to signatures, we recall some
standard results concerning real forms (see \cite{abv}, Chapter Two).
Recall
\[
\Gmin = \{ x\in G: \, \dtheta(x)=x^{-1} \}.
\]
Note that $G$ acts on $\Gmin$ by 
\[ 
g\star x=gx\dtheta(g^{-1})
\]
It is routine to check that if $\tau = {\rm Ad}_x   \dtheta$ is an
involution, then 
\[ 
{\rm Ad}_g  {\rm Ad}_x \dtheta   {\rm Ad}_{g^{-1}}=
{\rm Ad}_{g\star x}  \dtheta
\]

\begin{lem}
\label{lem_Hconj}
{\it
If $x\in \Gmin,$ there exists $g\in G$ such that $g\star x=t\in T^{\gamma_d}$ 
is of order $2.$
}
\end{lem}

\noindent
{\bf Proof.} Since $x\in \Gmin,$ ${\rm Ad}_x  \dtheta$ is an involution.
By conjugating  in $G,$ we may assume ${\rm Ad}_x   \dtheta$ and $\theta$ 
commute (see \cite{hel}, Theorem III.7.1 and following remark).
It follows that $x\in K.$ By \cite{karo:homog-compact}, 
there exists $u\in K$ such that $u\star x\in T^{\gamma_d},$
so $u\star x\in T^{\tau_d}$ since $\theta$ acts trivially on $T$.
But $u\star x\in \Gmin,$ so $u\star x= (u\star x)^{-1},$ and
hence $u\star x$ is of order two.
\qed

\begin{lem}
\label{lem_Gconjform}
Any $\fl \in \Lagr(\fg, d)$
is $G$-conjugate to a real form $\fl_{d,\sigma}$ for
some $\daut$-signature $\sigma.$
\end{lem}

\noindent
{\bf Proof.} We know any real form in $\Lagr(\fg, d)$
is of of the form 
${\rm Ad}_g  \dtheta$ with $g\in \Gmin,$ so by the previous lemma, by
$G$-conjugation it
can be put in the form ${\rm Ad}_t  \dtheta$ for some $t\in T^{\gamma_d}$
of order $2.$ Since $t$ is of order $2,$ the eigenvalue $\sigma_t(\alpha)$
of $t$ on ${\fg}_\alpha$ is $\pm 1.$ It is easy to check that $\sigma_t$
is a signature, and since $t\in T^{\gamma_d},$ it is a $\daut$-signature.
Hence our real form is conjugate to $\sigma_t  \dtheta.$
\qed

\subsection{The $G$-orbit of $\fls$}
\label{subsec_normfls}

\bigskip
Let  $\sigma$ be an extended $\daut$-signature $\sigma$
 with $\supp(\sigma)=[S_\sigma]$. 
We will use $L_{d, \sigma}, M_{\sigma}, P_\sigma$
and $N_\sigma$ to denote the connected subgroups of $G$ with
Lie algebras $\fl_{d, \sigma}, \fm_\sigma, \fp_\sigma$ and
$\fn_\sigma$ respectively. Recall also that
$\fz_\sigma = \fz_{S_\sigma}$.

\begin{lem}
\label{lem_meps}
\[ 
\dim_{\Bbb R} G\cdot \fls = \dim_{\Bbb C} \fg - \dim_{\Bbb C} \fz_\sigma.  
\]
\end{lem}

\noindent {\bf Proof. } 
This follows from Lemma \ref{lem_lagrnormalizer}.
\qed


\bigskip



\begin{lem}
\label{lem_lu-put}
Let $\sigma$ be an extended $\daut$-signature with $d$ trivial.
 Then $G\cdot \fls = K\cdot A\cdot\fls.$
\end{lem}

\noindent
{\bf Proof.} Since $K$ acts transitively on $G/P_\sigma$ and $P_\sigma$
has Levi decomposition $P_\sigma=M_\sigma N_\sigma,$ we can
write $g=kmn,$ $k\in K, m\in M_\sigma, n\in N_\sigma.$
For any real reductive group $G$ and the fixed point subgroup $G_0$  of
an involution, there is a Cartan decomposition 
$G=K\tilde{A}G_0$ 
 where $\tilde{A}$ is chosen so that its Lie algebra $\tilde{\fa}$ has maximal
intersection with ${\fg}^{-\sigma,-\theta}$, the subspace of $\fg$ on
which $\sigma$ and $\theta$ act as $-1$
 (see \cite{ross:symspaces}, Theorem 10). When $d$ is trivial,
any real form $G_0$ in $\Lagr (\fg,d)$ contains a Cartan subalgebra
of $\fk$, so up to $K$-conjugacy we can choose $\tilde{\fa}=
\fa=i\ft$ so we can take $\tilde{A}=A$, the Iwasawa factor.
By the Cartan decomposition applied to the group $M_\sigma,$
we can write $m=k_max,$ with $k_m\in M_\sigma\cap K, a\in A, x\in 
M_\sigma^\theta.$ Thus, we can write $g=k_1au,$ with
$k_1\in K, a\in A, u\in L_{d,\sigma}.$ 
\qed

\section{$\overline{\Lagr(\fg, d)}$ as the real part of 
the De Concini-Procesi compactification $Z_d$ of $G$}
\label{sec_dp}

\bigskip
In this section, we identify the variety $\overline{\Lagr(\fg,d)}$
with the real points of a De Concini-Procesi
compactification $\Zdaut$ of the group $G.$ Since 
$\Zdaut$ is known to be smooth, it follows that $\overline{\Lagr(\fg,d)}$
is a manifold. We also show that $\overline{\Lagr(\fg,d)}$ is
connected and 
determine the $G$-orbits
in $\overline{\Lagr(\fg,d)}.$

\subsection{The complexification of $\fg$}
\label{subsec_cmplx}

\bigskip
Regard $\fg$ as a real Lie algebra and denote its complex structure
by $J_0 \in {\rm End}_{\Bbb R}(\fg)$. 
We may identify its complexification  ${\fg}_{\Bbb C}$ with
$ (\fg \oplus \fg, J_0 \oplus J_0)$ via the map
\[
{\fg}_{\Bbb C} \, \lrw \, (\fg \oplus \fg, J_0 \oplus J_0): \, \,
x + iy \, \Map \, (x + J_0 y, \, \theta(x) + J_0 \theta(y)), 
\hspace{.1in} x, y \in \fg.
\]
Under this identification, the complex conjugation
operator $\tau$ on $\fg_{\Bbb C}$ becomes 
\[
\tau(X,Y)=(\theta(Y),\theta(X)),
\]
with its set of real points realized as
\[
(\fg\oplus \fg)^{\tau}=\{ (X, \theta(X)):  X\in \fg \}.
\]
If $\fr\subset \fg$ is a real subalgebra, then ${\fr}_{\Bbb C}=\fr + i\fr$
is regarded as 
 a complex subalgebra of $\fg\oplus \fg.$ For example, ${\fk}_{\Bbb C}$
is the diagonal subalgebra ${\fg}_\Delta = \{ (X,X):  X\in \fg \}$
and $(\ft+\fnp)_{\Bbb C}={\fh}_\Delta + {\fn}_1 + {\fnm}_2,$
where
for a Lie subalgebra $\fr$ of $\fg,$ 
\begin{equation}
\label{eq_ggnot}
{\fr}_\Delta = \{ (X,X): \,  X\in \fr \}, \ \ {\fr}_1=
\{(x, 0): x \in \fr\}, \, \, 
\fr_2 = \{(0, x): x \in \fr\}.
\end{equation}
The proof of the following lemma is straightforward.

\begin{lem}
\label{lem_lcplx}
For an extended $d$-signature $\sigma$, the complexification
$\flsC$ of $\fls$ is 
\[
{\flsC}= \{(X, \, a_\sigma\dgamma(X)): X\in \fm_\sigma \} \oplus
{\fn_{\sigma}}_1 \oplus {\fnsigmaminus}_2
\]
\end{lem}

\bigskip
Recall that $\ll \, , \, \gg$ is the Killing form of $\fg$.
Consider the symmetric form $I$ on $\fg\oplus \fg$ given by
\[
I((x_1,x_2),(y_1,y_2))=\ll x_1,y_1\gg- \ll x_2,y_2\gg.
\]
Then $\fl\subset \fg$ is a real Lagrangian subalgebra of $\fg$
with respect to the imaginary part of the Killing form if and only
if ${\fl}_{\Bbb C}\subset \fg\oplus \fg$ is a complex Lagrangian subalgebra
with respect to $I.$ 

If we denote  by ${\cal L}_{\Bbb C}$
the set of all complex Lagrangian subalgebras
of $\fg \oplus \fg$ with respect to $I$, then we have the injective
map
\[
\Lagr \lrw {\cal L}_{\Bbb C}: \, \fl \Map \fl_{\Bbb C}.
\]
With respect to the Adjoint action of $G$ on $\Lagr$, we have
\[
\left({\rm Ad}_g \fl \right)_{\Bbb C} \, = \, {\rm Ad}_{(g,\theta(g))}
(\fl_{\Bbb C}).
\]
 
On the group level, we have the analogous identification
$G_{\Bbb C}\cong G\times G$. We lift $\tau$ to an involution
also denoted $\tau$ of $G\times G$. 
 In this context,
$G$ (as the set of real points) is identified with the fixed point set of
$\tau$ as
\[ 
\{ (g, \, \theta(g)): \,  g\in G \}
\]
 Let $G_{\Delta, d}=\{(x,\dgamma(x)): \,  x\in
G\}$.
Then $(G\times G)/G_{\Delta, d}$ is
an example of a complex symmetric space, and 
De Concini and Procesi \cite{dp:compactification}
have exhibited a particular smooth compactification $\Zdaut$ 
of $(G\times G)/G_{\Delta,d}.$

\subsection{The De Concini-Procesi compactification $Z_d$}
\label{subsec_dpcmpt}

\bigskip
Note that $G\times G$ acts on the Grassmannian of
$n$-dimensional complex subspaces of  $\fg\oplus \fg$
through the Adjoint action, where $n = \dim_{\Bbb C} \fg$.
Consider the $\dgamma$-diagonal subalgebra 
\[
{\fg}_{\Delta,d} = \{ (X,\dgamma(X)):  X\in \fg \}
\]
of $\fg\oplus \fg$ and 
the orbit  $(G\times G)\cdot
{\fg}_{\Delta,d}$ inside the Grassmannian. The stabilizer subgroup
of $G \times G$ at $\fg_{\Delta, d}$ is $G_{\Delta, d}$, so
$(G\times G)\cdot
{\fg}_{\Delta,d} \cong (G \times G) /G_{\Delta, d}$.
By definition, the De Concini-Procesi variety is the 
closure (with respect to the Zariski or the classical topology)
of $(G\times G)\cdot
{\fg}_{\Delta,d}$ in the Grassmannian. It will be denoted by 
$Z_d$ and it is called the De Concini-Procesi compactification
(of $(G \times G) /G_{\Delta, d}$). It is a 
smooth complex manifold of complex dimension
$n$ (see \cite{dp:compactification} for more details). 
Since the variety of complex Lagrangian subalgebras is $G\times G$ stable,
it follows that every element in $Z_d$ is a complex Lagrangian
 subalgebra of $\fg \oplus \fg$
of dimension $n$.

\bigskip
It is known  \cite{dp:compactification}
that $G\times G$ has finitely many orbits on
$\Zdaut.$ We describe the orbits.
Recall that 
$S(\Sigma_+)=\{ \alpha_1,\dots, \alpha_l \}$ is the set of all simple roots.
Let
$\eta: S(\Sigma_+)  \to \{ 0,1 \}$ be any map. Regarding $\eta$
as 
an extended signature
for the trivial involution, we have
the parabolic subalgebra
\[
\fp_\eta=\fm_\eta+\fn_\eta
\]
and $\fnetaminus = \theta(\fn_\eta)$ of $\fg$.
Consider the subalgebra
\[
\fgdeta = \{(X,\dgamma(X)): \, X\in \fmeta \} \oplus {\fn_\eta}_1\oplus
\dgamma {\fnetaminus}_2 .
\]
Note that when $\eta$ is constant on
$\daut$-orbits and is regarded as an extended $d$-signature, 
we have  $\fgdeta ={\fl}_{d,\eta, {\Bbb C}}.$

\begin{thm}
\label{thm_GGorbits}
 \cite{dp:compactification} {\it
Every point $\fr \in \Zdaut$ is in a $G\times G$ orbit of $\fgdeta$ for some
$\eta.$
}
\end{thm}

\bigskip
We say that a complex
subalgebra $\fr$ of $\fg\oplus \fg$ has a real structure if
it is the complexification of a real subalgebra of $\fg$ under the
identification $\fg_{\Bbb C} \cong \fg \oplus \fg$.
This is equivalent to the condition that $\tau(\fr)=\fr,$ and
in this case,
\[
\fr=({\fr}^\tau)_{\Bbb C},
\]
where $\fr^\tau \subset \fg \oplus \fg$, the fixed point set of $\tau$ in $\fr$,
 is identified with its image in $\fg$ under the projection
$\fg \oplus \fg \rightarrow \fg: (x, y) \mapsto y$.

\begin{nota}
\label{nota_ZdR}
{\em
We will denote the set of all Lie algebras in
$\Zdaut$ with a real structure by ${\ZdR}.$
}
\end{nota}

Note that $\fg_{\Delta,d} \in Z_{d, {\Bbb R}}$. In fact, 
\[
\fg_{\Delta, d} \, = \, (\fl_{d, \sigma_1})_{\Bbb C},
\]
where 
$\sigma_1(\alpha) = 1$
for all $\alpha$.

In fact, $\fgdeta$ is in $Z_{d, {\Bbb R}}$ if and only if $\eta$ is
constant on $d$-orbits.

 Since $\tau$ preserves $\fg_{\Delta,d}$,
$\tau$ preserves the open  subset $(G \times G) \cdot \fg_{\Delta,d} \subset
Z_d$. 
Since $\tau$ is continuous,
it follows that $\tau$ preserves $\Zdaut.$ Thus, $\ZdR$ is the
set of real points of a complex compact manifold, so $\ZdR$ is a 
 compact manifold.

\subsection{$G$-orbits on $\ZdR$}
\label{subsec_Gorbits}
Recall that for every Lagrangian subalgebra $\fl\subset \fg$,
\begin{equation}
\label{cplx_equiv}
({\rm Ad}_g(\fl))_{\Bbb C}=({\rm Ad}_g, {\rm Ad}_{\theta(g)})({\fl}_{\Bbb C}),
\hspace{.2in} \forall g \in G.
\end{equation}


\begin{prop}
\label{prop_conjtocl}
Every $\fr\in \ZdR$ is $G$-conjugate to
${\flsC}$ for some extended $\daut$-signature $\sigma.$
\end{prop}

\noindent
{\bf Proof.} Let $\fr=(g_1,g_2)\cdot \fgdeta$ for some $\eta,$
 so
 \[
\fr=\{ ({\rm Ad}_{g_1}(y+z_1),
{\rm Ad}_{g_2} \dgamma(y+z_2) ): \,  
y\in \fmeta, z_1\in \fn_\eta, z_2\in \fnetaminus \}.
\]
Since $\fr$ has a real structure, $\tau (\fr)=\fr,$ so
\[
({\rm Ad}_{\theta(g_2)} \dtheta(y+z_2), \, 
{\rm Ad}_{\theta(g_1)}\theta(y+z_1))
\]
is in $\fr,$ so that ${\rm Ad}_{\theta(g_2)}\dtheta(y+z_2)=
{\rm Ad}_{g_1}(u+v_1)$
for  some $u\in \fmeta$ and 
$v_1\in \fn_\eta.$ But 
\[
\fpeta=\{ \theta(y+z_2): \,  y\in \fmeta, z_2\in
\fnetaminus \},
\]
so ${\rm Ad}_{g_1^{-1}\theta(g_2)}\dgamma(\fpeta)\subset \fpeta.$
Since $\dgamma(\fpeta)$ is $G$-conjugate to $\fpeta,$ it follows that
 $\dgamma(\fpeta)=\fpeta.$ 
Since $\Peta$ is the normalizer of $\fpeta,$ it follows that 
$g_1^{-1}\theta(g_2)\in \Peta,$ so
$g_2=\theta(g_1p),$
for some $p\in \Peta.$  Thus, 
\[
\fr=\{ ({\rm Ad}_{g_1}(y+z_1), \, 
{\rm Ad}_{\theta(g_1)\theta(p)}\dgamma(y+z_2) ): \,  
 y\in \fmeta, z_1\in \fn_\eta, z_2\in \fnetaminus \}.
\]
Thus, up to $G$-conjugacy, 
\[
\fr=\{ (y+z_1),
{\rm Ad}_{\theta(p)}\dgamma(y+z_2) ): \, 
  y\in \fmeta, z_1\in \fn_\eta, z_2\in \fnetaminus \}
\]
and $\fmeta,$ $\fn_\eta$ and $\fnetaminus$ are $\dgamma$-stable.

\bigskip
We write $\theta(p)=lu$ with $l\in \Meta,$ $u\in \Netaminus $. Since
\[
\{ u\cdot (y+z_2): \, y\in \fmeta, z_2\in \fnetaminus \} 
  = \{ (y+w_2): \,   y\in \fmeta, w_2\in \fnetaminus \} 
\]
 it follows that
\[
\fr=\{ ((y+z_1),
{\rm Ad}_l\dgamma(y+z_2) ): \, 
y\in \fmeta, z_1\in \fn_\eta, z_2\in \fnetaminus \} .
\]
We use again the assumption that $\fr$ has a real structure and the facts that
$\theta(\fmeta)=\fmeta$, $\theta(\fn_\eta)=\fnetaminus$, $M_\eta$
preserves the decompositions $\fpeta=\fmeta+\fn_\eta$ and
$\theta(\fpeta)=\fmeta+\fnetaminus.$
Since 
\[
\tau(y+z_1,{\rm Ad}_l\dgamma (y+z_2))=
({\rm Ad}_{\theta(l)}\dgamma(\theta(y)+\theta(z_2)), \,
\theta(y)+\theta(z_1)),
\]
we see that 
\[
\{{\rm Ad}_{\theta(l)}\dgamma y, \, y): \, y\in \fmeta \} =
\{ (y,{\rm Ad}_l\dgamma(y)): \, y\in \fmeta \}
 = \{ (\dgamma ({\rm Ad}_{l^{-1}} y), \, y): \, y\in \fmeta \}.
\]
Hence, ${\rm Ad}_{\theta(l)} \dgamma = \dgamma {\rm Ad}_{l^{-1}},$ and
it follows that $\dtheta(l)=l^{-1}.$

\bigskip
Now, by Lemma  
\ref{lem_Hconj}, there exists $v\in M_\eta$ such that $v\star l = t
\in T^{\gamma_d}$ of order $2.$ But it is easy to check that 
\[ 
(\theta(v),v)\cdot (1,l)\cdot \fgdeta = (1,v\star l) \fgdeta.
\]
Hence, after acting by an element of $M_\eta,$ 
we may assume that $\fr=(1,t)\cdot \fgdeta$ and that $t\in T^{\gamma_d}$
 is an element of
order $2.$

\bigskip
As before, let $\sigma_t (\alpha)$ be the eigenvalue of $t$ on the root space 
${\fg}_\alpha.$ Then $\sigma_t$ is a
  $\daut$-signature
and we can define a new 
extended $\daut$-signature ${\sigma}^\prime$ by 
\[
{\sigma}^\prime (\alpha)=\eta(\alpha)\sigma_t(\alpha)
\]
Then $(1,t)\cdot \fgdeta =\fl_{{d,{\sigma}^\prime}, {\Bbb C}},$
using Lemma \ref{lem_lcplx}, which completes the proof of
the proposition.
\qed

\subsection{Geometry and topology of the closure $\overline{\Lagr(\fg, d)}$}
\label{subsec_geoclose}

\begin{thm}
\label{thm_zrconn}
${\ZdR}$ is connected.
\end{thm}

\noindent
{\bf Proof. } Since $G$ is connected, Proposition \ref{prop_conjtocl}
implies that
it suffices to find a path from ${\flsC}$ to the
solvable Lie algebra 
${\flsnullC}$, where $\sigma_0(\alpha)=0$ for all $\alpha \in \Sigma$.
Note that
\[
{\flsnullC}=\{ (H,\dgamma(H)): \, H\in \fh \} \oplus {\fn}_1\oplus
{\fnm}_2.
\]
Let $H\in \fh$ have the property that $\alpha(H)>0$ for all $\alpha
\in {\Sigma}_+.$ 
If $X\in {\fm}_\sigma \cap \fg_\alpha, \alpha \in {\Sigma}_+$ then 
\[
{\lim}_{t\to +\infty} ({\rm Ad}_{\exp(tH)}, \, 
{\rm Ad}_{\theta(\exp(tH)}) {\Bbb C }(X,\dgamma(X))
={\Bbb C} (X,0),
\]
and if $X\in \fm_\sigma \cap \fg_\alpha, \alpha \in -{\Sigma}_+,$
\[
{\lim}_{t\to +\infty} ({\rm Ad}_{\exp(tH)},
\, {\rm Ad}_{\theta(\exp(tH)}) {\Bbb C} (X,\dgamma(X))
={\Bbb C} (0,\dgamma(X)).
\]
Since 
\[
{\flsC} = \{ (X,\dgamma_\sigma(X)): \,  X\in {\fm}_{\sigma} \}
 \oplus {\fn_\sigma}_1
\oplus {\fnsigmaminus}_2,
\]
it follows that 
\[
{\lim}_{t\to +\infty} ({\rm Ad}_{\exp(tH),\theta(\exp(tH))}){\flsC} 
= {\flsnullC}.
\]
\qed

\begin{rem}
\label{rem_closedorbit}
{\em
This theorem can also be proved by observing that $\ZdR$ has a 
unique closed $G$-orbit $G\cdot \flsnullC$. The Lie algebra
$\fl_{d,\sigma_0}={\fh}^{\tau_d} + \fn$. When $d$ is non-trivial,
and $\sigma$ is a $d$-signature, the curve $\exp({\rm ad}_{tH}) \
\fl_{d,\sigma}$
provides a class of examples when $T: \Lagr \rightarrow \Lagr$ 
is not continuous (see
Remark \ref{rem_notcont}).
}
\end{rem}

\begin{nota}
\label{nota_Zdeta}
{\em
We will use $\Zdeta$ to denote the $G\times G$-orbit
 through $\fgdeta$. 
We let $\eta_1$ be the extended $\daut$-signature such that
$\eta_1(\alpha_i)=1,$ all $\alpha_i \in S(\Sigma_+).$ Then ${\fg}_{d,{\eta}_1}
= {\fg}_{\Delta,d},$ and ${\Zdetao}$ is the unique open $G\times G$
orbit in $\Zdaut.$
}
\end{nota}

\begin{thm}
\label{thm_lageqzdr}
\[
\overline{\Lagr (\fg,\daut)}\cong {\ZdR}
\]
under the complexification map $\fl\to {\fl}_{\Bbb C}.$
In particular, $\overline{\Lagr (\fg,\daut)}$ is a smooth manifold.
\end{thm}

\noindent
{\bf Proof. } 
By Proposition \ref{prop_conjtocl}, we know $G$ has finitely many orbits
on ${\ZdR},$ and the orbits are 
given by extended $\daut$-signatures. The open orbits
are given by the orbits through ${\flsC},$ where $\sigma$ is a
$\daut$-signature. Indeed, in the proof of 
Proposition \ref{prop_conjtocl}, we
showed that ${\Zdeta}\cap {\ZdR}$ is a finite disjoint
union of $G$-orbits $G\cdot {\flsC}$ with $|\sigma(\alpha)|=
\eta(\alpha)$ for every root $\alpha.$
 Moreover, each of these $G$-orbits has the same dimension by Lemma
 \ref{lem_meps}.
 It follows that the orbits $G\cdot {\flsC}$
are the connected components of ${\Zdeta} \cap {\ZdR}$
for $\eta = |\sigma|$ 
and also that the $G\cdot {\flsC}$ are locally closed. Since
${\Zdetao}$ is open, the orbits $G\cdot {\flsC}$ are
open when $\sigma$ is a $\daut$-signature, and by the dimension
statement, none of the other orbits are open since ${\ZdR}$ is
connected. Moreover, it follows from the fact that ${\ZdR}$
is a finite union of locally closed orbits that the union of the
open orbits is dense.

\bigskip
Now it suffices to prove that $\Lagr(\fg,d)$ surjects onto the open orbits
of ${\ZdR}.$ By Lemma 
 \ref{lem_Gconjform}, we know that every real form in
$\Lagr(\fg,\daut)$ is ${\rm Ad}_g{\fls},$ for some
 $\daut$-signature $\sigma.$ 
It follows from  (\ref{cplx_equiv}) and
the above description of open orbits on ${\ZdR}$
that $\Lagr(\fg,\daut)$
  maps onto the union of the open orbits of ${\ZdR}.$
\qed

\begin{lem}
\label{lem_zclosed}
The Zariski closure of $\Lagr(\fg,d)$ coincides with its closure
in the classical topology.
\end{lem}

\noindent{\bf Proof. } We know
\[
\Lagr(\fg,d)=\cup_{\sigma} G\cdot {\fls}
= {\Zdetao} \cap {\ZdR},
\]
 where the union is over all $d$-signatures and $\Zdetao$ is
the open $G\times G$ orbit on $\Zdaut.$ Thus, $\Lagr(\fg,d)$ is the
real points of $\Zdetao.$ But the Zariski closure of the real points
is contained in the real points of the Zariski closure, so the Zariski closure
of $\Lagr(\fg,d)$ is contained in ${\ZdR}=\overline{\Lagr(\fg,d)}.$
Since the classical closure of $\Lagr(\fg,d)$ is
contained in the Zariski closure, it follows that they coincide.\qed

\subsection{Open orbits in $\overline{\Lagr(\fg,d)}$}
\label{subsec_openorbits}

\bigskip
In this subsection we identify the open orbits in 
$\ZdR$ with symmetric spaces. 

\begin{prop}
\label{prop_norm}
Let $\tau$ be a real form of a semisimple Lie algebra $\fg$, and also
denote its lifting to the adjoint group $G$ by $\tau$.
Then $$G^\tau=N_G({\fg}^\tau).$$ (see \cite{dp:compactification} 
for the holomorphic
version of this fact. The proof is essentially the same).
\end{prop}

\begin{cor}
\label{cor_symm}
{\it 
For a $\daut$-signature $\sigma,$ the open orbit $G\cdot {\flsC}$
 is the semisimple symmetric space
 $G/G^{\tau_{d,\sigma}}$.
}
\end{cor}

\noindent {\bf Proof.}
The above proposition implies that the stabilizer 
$N_G(\fls)=G^{\tau_{d,\sigma}}$.
\qed

\subsection{Another description for $\overline{\Lagr(\fg, {\rm id})}$}
\label{sec_zeLagr}

The set $\overline{\Lagr(\fg, {\rm id})}$ 
has been most important for applications. In this
section, we give another description of it.

\bigskip
When the diagram automorphism $\daut$ is trivial, we will refer to the
corresponding real De Concini-Procesi compactification as $Z_{\Bbb R}$ 
instead
of ${\ZdR}.$ By Theorem \ref{thm_lageqzdr}, $Z_{\Bbb R}=
\overline{\Lagr(\fg,  {\rm id})}$.
It will follow from the description of irreducible components in
Section \ref{subsec_irrcomp} that $Z_{\Bbb R}$ is
the unique irreducible component of $\Lagr$ containing $\fk.$
We let
\[
\zeLagr = \{ \fl \in \Lagr: \, \rank (\fk\cap \fl) = \rank (\fk) \}.
\]
It is the set of Lagrangian subalgebras of $\Lagr$ containing the
Lie algebra of a maximal torus of $\fk.$

\bigskip
\begin{prop}
\label{prop_isinzr}
$\zeLagr = Z_{\Bbb R}.$
\end{prop}

\noindent
{\bf Proof. }
Write $\fl_\sigma$ for $\fl_{d, \sigma}$ for $d$ trivial.
First assume $\fl={\rm Ad}_g\flds$ lies in $Z_{\Bbb R}.$
By Lemma \ref{lem_lu-put}, we can write $\fl=
{\rm Ad}_k {\rm Ad}_a{\fl}_\sigma,$ for $k\in K,$ $a\in A$.
 But ${\fl}_\sigma$ contains $\ft,$
so ${\rm Ad}_k {\rm Ad}_a{\fl}_\sigma$ contains ${\rm Ad}_k(\ft),$ since
$A$ acts trivially on $\ft.$ 
Thus, $\fl\in \zeLagr.$

\bigskip
Now assume a Lagrangian subalgebra $\fl$
 contains the Lie algebra of a maximal torus of
$K.$ By \cite{karo:homog-compact}, we know 
$\fl={\rm Ad}_k(\fm_{{\scriptscriptstyle S}, 1}^{\tau} 
\oplus V \oplus \fnS),$ for some $(S,V,\tau).$
By the assumption on $\fl,$ we may assume that 
$\fm_{{\scriptscriptstyle S}, 1}^{\tau} \oplus V$
contains $\ft.$ Then $V=\ft \cap \fzS$ and $
\fm_{{\scriptscriptstyle S}, 1}^{\tau} $ contains
a Cartan subalgebra of $\fmSss \cap \fk.$ But it is easy to show that if
$\tau$ does not have trivial diagram automorphism, then 
$\fm_{{\scriptscriptstyle S}, 1}^{\tau}$
does not contain a Cartan subalgebra of $\fmSss \cap \fk.$ It follows
easily that $\fl \in Z_{\Bbb R}.$
\qed

\bigskip
We remark that it follows that $G$ acts on $\zeLagr$, a fact that is
not clear from the definition of $\zeLagr.$

\bigskip
\begin{cor}
\label{cor_ZRmodelpoints}
All points of $Z_{\Bbb R}$ are model points
\end{cor}

\noindent
{\bf Proof. } 
This follows from Proposition \ref{prop_modelpt} and the observation
that if $\fl(S,V,\tau)$ contains $\ft$, then $V\subset \ft$.
\qed

\begin{rem}
\label{rem_poistrinZR}
{\em
It follows from Corollary \ref{cor_ZRmodelpoints}
 that many familiar Poisson structures
are contained in $Z_{\Bbb R}$ as $G$ or $K$ orbits with the Poisson
structures being the restriction of the Poisson
structure $\Pi$ on ${\cal L}$ defined in
Section \ref{sec_Lagr-general}.
For example, 
we can identify $G\cdot \fk\cong G/K$, and the Poisson
structure induced by $\Pi$  on $G/K \cong AN$
is  the negative of the Poisson structure $\pi_{AN}$ that makes
$AN$ into the dual  Poisson Lie group of $K$.
More generally, by looking at $G$-orbits in $\Lagr(\fg, d)$,
we obtain in this manner a Poisson structure on $G/G_0$ for every real
form $G_0$ of $G$. The Poisson manifolds arising from $K$-orbits
in $Z_{\Bbb R}$ are studied in more detail in Section 
\ref{sec_poi-on-lagr}.
}
\end{rem}

\begin{rem}
\label{rem_modelspois}
{\em
Not all points in $Z_{d,\Bbb R}$ are model points when $d$ is not
trivial. The criterion for ${\fl}_{d,\sigma}$ to be a model point is that
if $\sigma(\alpha)=0$, then $d(\alpha)=\alpha$. 
}
\end{rem}

\bigskip
In \cite{e-l:harm}, we introduced certain $K$-invariant metrics
$g_\lambda$ on $T^*(K/T)$ for $\lambda \in \fa_r$, the set of elements
in $\fa$ whose centralizer in $K$ is $T$. These metrics are
important for showing that an operator $S$ introduced by Kostant
is a limit of some Hodge Laplacians $S_\lambda.$ 
The existence of this family
simplifies the proof of Kostant's basic result that ${\rm Ker}(S)$
is isomorphic to $H^*(K/T).$ We remark that the metrics $g_\lambda$
can be understood in terms of the restriction of a Riemannian metric on
the Riemannian symmetric space $G/K.$ 
Since $Z_{\Bbb R}$ is a 
compactification of $G/K$ with closed orbit the flag manifold
$G/B,$ this observation provides evidence that embedding the
Bruhat-Poisson structure on $G/B$ into the manifold $Z_{\Bbb R}$ is
useful in Poisson geometry.

\bigskip
We give the construction of this metric. We can identify the tangent
space of $G/K$ at $gK$ with ${\rm Ad}_g(i\fk).$ The Killing form is
positive definite at ${\rm Ad}_g(i\fk),$ and we let $s$ be the metric
on $G/K$ given by taking the square root of the Killing form metric
on ${\rm Ad}_g(i\fk).$

\bigskip
Let $H_{\lambda}\in \fa $ be such that $\lambda(H)=(H_\lambda,H)$ and
let $a_\lambda =\exp(H_\lambda).$ Then the $K$-orbit through
$a_\lambda K\in G/K$ can be identified with $K/T.$ 
If we restrict the above metric $s$ to a metric $s_{\lambda}$
on  $K\cdot a_\lambda K\subset G/K,$
and use $s_{\lambda}$ to identify the cotangent bundle with the
tangent bundle, then one can show by easy calculations that
we obtain the metric $g_\lambda$ from \cite{e-l:harm}.

\section{Geometry of $\overline{\Lagr (S,\eps,d)}$}
\label{sec_lsedclos}

In this section, we combine results from Section \ref{sec_oshima} with results
from Section \ref{sec_dp} to study the closures $\overline{\Lagr (S,\eps,d)}.$

\subsection{Smoothness of $\overline{\Lagr(S,\eps,d)}$}
\label{subsec_smoothclos}

\begin{thm}
\label{thm_smoothcomp}
Each $\overline{\Lagr (S,\eps,d)}$ is a smooth connected 
submanifold of the Grassmannian $\Grg$
 of dimension
$\dim(\fk)+{{z(z-3)}\over 2}.$ It fibers over $G/\PS$ with  the fiber
being the product of $\LagrfzSeps$ with $\overline{\Lagr(\fmSss,\daut)}$,
the real points of a De Concini-Procesi variety. 
\end{thm}

\noindent
{\bf Proof. } 
Recall from the proof of Proposition \ref{prop_lsedfiber} that
\[
\Lagr_{\fpS}(S,\eps,d)
\cong \LagrfzSeps \times \Lagr(\fmSss,\daut)
\]
Thus,
\[
\overline{\Lagr_{\fpS}(S,\eps,d)}\cong \LagrfzSeps \times
\overline{\Lagr(\fmSss,\daut)}
\]
because $\LagrfzSeps$ is already closed.
Once we identify
\[
K\times_{K\cap P_S}\overline{\Lagr_{\fpS}(S,\eps,d)}\cong 
\overline{\Lagr (S,\eps,d)}
\]
the theorem will follow from Theorem \ref{thm_lageqzdr},
 Theorem \ref{thm_zrconn}, and Proposition \ref{prop_indlagr}.

So we consider the map 
\[
m:K\times_{K\times P_S}\overline{\Lagr_{\fpS}(S,\eps,d)}\to
\overline{\Lagr (S,\eps,d)}
\]
given by $m(k,\fl)=Ad_k \fl$.
It is easy to see that Karolinsky's Theorem \ref{thm_karolinsky}
implies that $m$ is onto. It suffices to check that $m$ is an
immersion, since $m$ is clearly smooth and proper. To show $m$ is
injective, suppose that for $i=1,2$,
$\fl(S_i,V_i,\tau_i)\in \overline{\Lagr_{\fpS}(S,\eps,d)}$
and $Ad_{k_1}\fl(S_1,V_1,\tau_1)=Ad_{k_2}\fl(S_2,V_2,\tau_2)$.
It follows as in the proof of Proposition \ref{prop_uniqueness}
 that $S_1=S_2$ and ${k_1}^{-1}k_2 \in K\cap P_{S_1}$.
Note that $\fnS \subset {\fn}_{\scriptstyle S_1}$, so $P_{S_1}\subset \PS$.
It follows easily that $m$ is injective, 
and the proof that the tangent map $m_*$ is injective
is similar to the proof of the same fact in Proposition
\ref{prop_lsedfiber}.

The dimension statement is clear from Proposition \ref{prop_lsedfiber}.
\qed

\subsection{Irreducible components}
\label{subsec_irrcomp}
\bigskip
In this subsection we determine the irreducible components of
$\Lagr.$

\begin{prop}
\label{prop_lsepsdirr}
$\overline{\Lagr(S,\eps,d)}$ is Zariski closed and irreducible.
\end{prop}

\noindent{\bf Proof.} 
Since $\overline{\Lagr(\fmSss,d)} \times \LagrfzSeps$ is Zariski
closed in $\Grg$ via the embedding $(\fl,V)\to \fl+V+\fnS$
 (Proposition \ref{prop_indlagr} and 
Lemma \ref{lem_zclosed}), it follows that
$G\times_{{\scriptscriptstyle{P_S}}}
 (\overline{\Lagr(\fmSss,d)} \times \LagrfzSeps)$ is
Zariski closed in $G\times_{{\scriptscriptstyle{P_S}}} 
{\rm Gr}(n,\fg)$. Moreover, the map
$m:G\times_{{\scriptscriptstyle{P_S}}}
 (\overline{\Lagr(\fmSss,d)} \times \LagrfzSeps)
\to {\rm Gr}(n,\fg)$ is projective, so its image is Zariski closed,
and irreducible since the domain is irreducible. Thus, the
proposition follows from Theorem \ref{thm_smoothcomp}.
\qed


\begin{dfn}
\label{dfn_iness}
\em{
Lagrangian data $(S,\eps,d)$ is  said to be
{\it inessential} if $S=S(\Sigma_+) - \{ \alpha_i \},$
$d=d^{'} |_S$ for some diagram automorphism $d^{'}$ of $S(\Sigma_+),$
and $\eps=1.$ Otherwise, $(S,\eps,d)$ is called essential.}
\end{dfn}

\begin{prop}
\label{prop_inesscon}
Lagrangian data $(S,\eps,d)$ is inessential if and only if
$\Lagr(S,\eps,d) \subset \partial \overline{\Lagr(S^{'},\eps^{'},d^{'})}$
for some Lagrangian data $(S^{'},\eps^{'},d^{'}).$
\end{prop}

\noindent{\bf Proof. } If $(S,\eps,d)$ is inessential, then we
claim $\Lagr(S,\eps,d)\subset \overline{\Lagr(S(\Sigma_+),1,d^{'})},$
 where $d^{'}|_S = d.$ Indeed, since $\dim (\fzS)=1$ and $\eps=1,$ the
Lagrangian subspace $V$ in $\fzS$ is $\fzS\cap \ft.$
It follows from Theorem \ref{thm_karolinsky} and  Lemma \ref{lem_Gconjform}
that each subalgebra in $\Lagr(S,\eps,d)$ is $G$-conjugate
to $\fm_{S,1}^{\tau_{d,\sigma}}\oplus \fzS \cap \ft \oplus \fnS$ for
some $\sigma.$ But this algebra coincides with ${\fl}_{d^{'},\sigma}.$ Hence,
$\Lagr(S,\eps,d) =\cup_{\sigma} G\cdot {\fl}_{d^{'},\sigma},$
so
\[
\Lagr(S,\eps,d)\subset Z_{d^{'},\Bbb R}=\overline{\Lagr(\fg,d^{'})}.
\]

Suppose that $\Lagr(S,\eps,d)\subset \partial 
\overline{\Lagr(S^{'},\eps^{'},d^{'})}.$
It follows that $S\subset S^{'}$ so $\dim(\fzS) > \dim(\fzSprime).$
Moreover, by Theorem \ref{thm_smoothcomp}, we have
\[
{{\dim(\fzS) (\dim(\fzS)-3)}\over 2} <
{{\dim(\fzSprime) (\dim(\fzSprime)-3)}\over 2}.
\]
It follows that $\dim (\fzS)=1$ and $\dim (\fzSprime)=0.$
Thus, $\Lagr(S^{'},\eps^{'},d^{'})=\Lagr(\fg,d^{'})$ consists of real forms.
But $\overline{\Lagr(\fg,d^{'})}=Z_{d^{'},\Bbb R},$ so every subalgebra
in $\Lagr(S,\eps,d)$ is $G$ conjugate to some $\fls$ by Proposition
\ref{prop_conjtocl}. Since $\fzS$ is one-dimensional, and ${\gamma}_{d^\prime}$
acts by permutations on $\fh$, it follows that ${\gamma}_{d^\prime}$ acts
trivially on $\fzS$, so the Lagrangian subalgebra of $\fzS$ associated
by Karolinsky's classification with $\fls$ is $\fzS \cap \ft$. Thus, 
$\fls \in \Lagr(S,1,d^{'}|_S),$ and the assertion follows.
\qed

\begin{cor}
\label{cor_irrcom}
\em{
\[ \Lagr = \cup_{\rm{essential } (S,\eps,d)} \overline{\Lagr(S,\eps,d)}
\]
is the decomposition of $\Lagr$ into irreducible components.}
\end{cor}

\noindent{\bf Proof.}
By Proposition \ref{prop_lsepsdirr}, each $\overline{\Lagr(S,\eps,d)}$
is irreducible. Thus, the irreducible components are the
$\overline{\Lagr(S,\eps,d)}$ not properly contained in any other
$\overline{\Lagr(S^{'},\eps^{'},d^{'})}.$ By the previous Proposition, these
correspond to essential data.
\qed

\begin{cor}
\label{cor_containsk}
$\overline{\Lagr(S(\Sigma_+), 1, {\rm id})} \cong 
\overline{\Lagr(\fg, {\rm id})} \cong \Lagr_0$ is the only irreducible
component of $\Lagr$ containing $\fk$.
\end{cor}

\noindent{\bf Proof.}
The Zariski closure of $G\cdot \fk$ is easily seen to be
 $\overline{\Lagr(S(\Sigma_+), 1, {\rm id})}$, which is not contained
in any other irreducible component by the previous Corollary.
\qed

\bigskip
Note also that $\Lagr$ itself is typically not smooth, because
different irreducible components can intersect. This does not happen
for ${\fsl}(2),$ but for ${\fsl}(3),$ the components
$\overline{\Lagr(S(\Sigma_+),1,{\rm id})}$ and 
$\overline{\Lagr(\emptyset,1,{\rm id})}$
intersect in the flag variety of $SL(3,\C).$

\section{The Poisson structure $\Pi$ on $\Lagr$}
\label{sec_poi-on-lagr}

In this section, we study some properties of the Poisson
structure $\Pi$ on $\Lagr$ defined in 
Section \ref{sec_general}. More specifically, we relate $\Pi$ to
the Bruhat Poisson structure and determine the 
$(K, \pi_K)$-homogeneous Poisson spaces defined by points in 
$\Lagr_0 \cong \overline{\Lagr(\fg, {\rm id})}$.

\subsection{The fibre projection $\overline{\Lagr(S, \epsilon, d)} 
\rightarrow G/\PS$ is Poisson}
\label{sec_fiber-poisson}

It is clear from the definition of $\Pi$ that every $G$-invariant
smooth submanifold of $\Lagr$ is a Poisson submanifold. Thus,
each $\overline{\Lagr(S, \epsilon, d)}$ is a Poisson submanifold.
On the other hand, equip $G/\PS$ with the Bruhat Poisson structure 
$\pi_\infty$,
which is the unique $(K,\pi_K)$-homogeneous Poisson structure on
$G/\PS$ that vanishes at the identity coset $e\PS$. 
Recall from Theorem \ref{thm_smoothcomp} that we have the fiber bundle
$\overline{\Lagr(S,\eps, d)}\to G/\PS$.

\begin{prop}
\label{prop_projection-posson}
The fiber projection $\phi$ from 
$\overline{\Lagr(S, \epsilon, d)}$ to $G/\PS$ is a Poisson map.
\end{prop}

\noindent
{\bf Proof.}
First, we observe that the projection $\phi$ is $G$-equivariant.
Indeed, we can identify 
$K\times_{K\cap P_S} \overline{\Lagr_{\fpS} (S,\eps,d)}$
with $G\times_{P_S} \overline{\Lagr_{\fpS} (S,\eps,d)}$ via the obvious
inclusion, and the map from
$G\times_{P_S} \overline{\Lagr_{\fpS} (S,\eps,d)}$
  to $\overline{\Lagr(S,\eps,d)}$ is given by
the Adjoint action $(g,\fl)\mapsto {\rm Ad}_g \fl$. Then the projection
to $G/\PS$ is given by $(g,\fl)\mapsto g\PS$, which is obviously
$G$-equivariant.

Recall that the Poisson structure on $\overline{\Lagr(S,\eps,d)}$
is induced by the element ${1\over 2}R\in \wedge^2 \fg$ given in Section
\ref{sec_Lagr-general}. Since $\phi$ is $G$-equivariant, it follows
that $\phi_* \Pi$ is given by the bi-vector field on $G/\PS$ induced
by ${1\over 2}R$, so we just have to check that ${1\over 2}R$ induces
the Bruhat Poisson structure on $G/\PS$. It follows from the definition
of the Drinfeld map that the Lagrangian subalgebra associated with
the point $e\PS$ by $\pi_\infty$ is $(\fk\cap \fpS)
\oplus \fnS$. By Theorem \ref{thm_M-to-O}, the Drinfeld map 
\[
\Dr:(G/\PS,\pi_\infty)
\to (K\cdot((\fk\cap \fpS)\oplus \fnS) ,\Pi)
\]
is a Poisson map. The normalizer of $(\fk\cap \fpS)\oplus \fnS$
in $K$ 
is $K\cap \PS$, and it follows that the Drinfeld map is a diffeomorphism,
so $\pi_\infty$ coincides with $\Pi$.
Since the Poisson structure $\Pi$ is induced by ${1\over 2}R$, the result
follows.
\qed

\subsection{$(K,\pi_K)$-homogeneous Poisson spaces determined by points
in $\Lagr_0$}
\label{sec_homog-poi-from-Lagrzero}

We now turn to the Poisson submanifold $(\Lagr_0, \Pi)$, where
$\Lagr_0 \cong \overline{\Lagr(\fg, {\rm id})}$ is the
unique irreducible 
component of $\Lagr$ that contains $\fk$.
We study the $(K, \pi_K)$-homogeneous Poisson spaces
determined by points in $\Lagr_0$ (see Definition \ref{dfn_standard-C}).

By Corollary \ref{cor_ZRmodelpoints}, every point in $\Lagr_0$ is a
model point. It follows from the discussion in Section 
\ref{sec_modelpt} that each $\fl \in \Lagr_0$  can determine
a number of $(K, \pi_K)$-homogeneous Poisson spaces, Indeed,
let $N_K(\fl)$ be the normalizer subgroup of $\fl$ in $K$.
Then for any subgroup $K_1$ of $K$ with the same Lie algebra 
$\fl \cap \fk$ as $N_K(\fl)$, the space
$K/K_1$ carries a unique Poisson structure $\pi$ such that the
covering map
\[
P: \, K/K_1 \lrw K/N_K(\fl) \cong K \cdot \fl \subset \Lagr_0: \, 
kK_1 \Map kN_K(\fl)
\]
is a Poisson map. The space $(K/K_1, \pi)$ is automatically
$(K, \pi_K)$-homogeneous, and the map $P$ is its Drinfeld map
(see Definition \ref{dfn_drinfi-lag}). 
Examples of $K_1$ are $K_1 = N_K(\fl)$ or $K_1$ is the
 connected component of the identity of $N_K(\fl)$.
 We can characterize these $(K, \pi_K)$-homogeneous Poisson spaces 
determined by points $\fl \in \Lagr_0$ as follows.

\begin{prop}
\label{prop_orbits-in-lagr-zero}
All $(K, \pi_K)$-homogeneous Poisson spaces $(K/K_1, \pi)$
determined by points in $\Lagr_0$ (see Definition \ref{dfn_standard-C})
have
the property that 
$K_1$ contains 
a maximal torus of $K$.
Conversely, all $(K,\pi_K)$-homogeneous Poisson spaces
with this property are determined by points in $\Lagr_0$.
\end{prop}

\noindent
{\bf Proof.} The first part of the proposition follows
from the definition of $\Lagr_0$. Now let
$(K/K_1, \pi)$ be any $(K, \pi_K)$-homogeneous Poisson
space such that $K_1$ contains
a maximal torus of $K$. Then  the Lie algebra $\fk_1$ of $K_1$
contains the
Lie algebra of a maximal torus of $K$.  Consider 
the Drinfeld map
\[
P: \, K/K_1 \lrw \Lagr.
\]
Let $\fl = P(eK_1) \in \Lagr$. Then by Drinfeld's Theorem
\ref{thm_drinfi-homog}, $\fk_1 = \fl \cap \fk$ 
and $K_1 \subset N_K(\fl)$. Thus $\fl \in \Lagr_0$ by the
definition of $\Lagr_0$, and $(K/K_1, \pi)$ is determined
by $\fl$.
\qed

The second part of  
Proposition  \ref{prop_orbits-in-lagr-zero}
can be rephrased as the following.

\begin{cor}
\label{cor_Kone}
Every $(K, \pi_K)$-homogeneous Poisson space
$(K/K_1, \pi)$, where $K_1$ is a closed subgroup
of $K$ containing a maximal torus of $K$,
is a Poisson submanifold of $(\Lagr_0, \Pi)$ up to 
a covering given by its Drinfeld map.
\end{cor}

\begin{rem}
\label{rem_Q}
{\em
Examples of $K_1$ in Proposition \ref{cor_Kone} are $K \cap Q$,
where $Q$ is a parabolic sungroup of $G$, so the corresponding
homogeneous space is a flag manifold $K/(K \cap Q) \cong G/Q$.
}
\end{rem}

\subsection{The normalizer subgroup of $\fl \in \Lagr_0$ in $K$}
\label{sec_normalizer}
 
We now study the normalizer subgroup 
$N_K(\fl)$ of an arbitrary $\fl \in \Lagr_0 $ in $K$ and determine 
when it is connected.
By 
Lemma \ref{lem_lu-put} and Proposition \ref{prop_conjtocl},
we can write
$\fl = {\rm Ad}_k {\rm Ad}_a \fl_{d,\sigma}$ for
some $k \in K, a \in A$ and extended signature $\sigma$
for $d = {\rm id}$, the trivial diagram automorphism. 
In what follows, we will write $\fl_\sigma = \fl_{{\rm id}, \sigma}$
and call an extended signature for $d = {\rm id}$ simply
an {\it extended signature}.
Write $a=\exp H
$ with $H\in \fa$ and further decompose
$H=H_1 + H_2$ with
$H_1\in \fa\cap\fmSigss$ and $H_2\in \fa\cap\fz_\sigma$.
 Then ${\rm Ad}_{\exp H}
\flstriv={\rm Ad}_{\exp H_1 }\flstriv$ since $H_2$
normalizes $\flstriv$.
Thus, we can assume $\fl= {\rm Ad}_{\exp H}\flstriv$
with $H\in \fa\cap \fmSigss$.
We will write $\fl_{H, \sigma} = 
{\rm Ad}_{\exp H} \fl_\sigma$.

\begin{lem}
\label{lem_lk}
For $\fl_{H, \sigma} = {\rm Ad}_{\exp H}\flstriv$,
where $\sigma$ is an extended 
signature and 
$H\in \fa\cap \fmSigss$,
\[
\fl_{H, \sigma} \cap \fk \, = \, \ft + \fn_\sigma + {\rm span}_{\Bbb R}
\{\Xa, \Ya: \, \sigma(\alpha) = 1, \alpha(H) = 0\}.
\]
\end{lem}

\noindent
{\bf Proof.} This follows from the fact that
\beqa
{\rm Ad}_{\exp H}\fl_\sigma \, = \, \ft + \fn_\sigma & + & 
{\rm span}_{\Bbb R} \{{\rm Ad}_{\exp H}\Xa, {\rm Ad}_{\exp H}\Ya: \sigma(\alpha) = 1\}\\& + & 
{\rm span}_{\Bbb R} \{i{\rm Ad}_{\exp H}\Xa, i{\rm Ad}_{\exp H}\Ya: \sigma(\alpha) = -1\}.
\eeqa
\qed

We now describe the normalizer subgroup of
$\fl_{H, \sigma}$ in $K$.

\begin{nota}
\label{nota_sigma}
{\em
For an extended signature $\sigma$ and $H \in \fa \cap \fm_{\sigma, 1}$,
let
$\Sigma_\sigma = \{\alpha \in \Sigma: \sigma(\alpha) = 1\}$. 
Let $W_\sigma$ be the subgroup of the Weyl group generated by the 
simple reflections corresponding to the simple roots in the
support of $\sigma$. Let
\[
W_{H, \sigma} \, = \, \{w \in W_\sigma: w \Sigma_\sigma = \Sigma_\sigma,
wH = H\} \subset W_\sigma \subset
W.
\]
Let 
\[
N^{'}(\fl_{H, \sigma}) \, = \, p^{-1}(W_{H, \sigma}),
\]
where $p: N_K(\ft) \rightarrow  W= N_K(\ft)/T$ is the projection
from the normalizer subgroup $N_K(\ft)$ of $\ft$ in
$K$ to the Weyl group. Finally, let
$K_{H, \sigma}$ be the connected subgroup of $K$ with
Lie algebra $\fl_{H, \sigma} \cap  \fk$.
}
\end{nota}

\begin{prop}
\label{prop_normalizer}
For an extended signature $\sigma$ and $H \in \fa \cap \fm_{\sigma, 1}$,
the normalizer subgroup $N_K(\fl_{H,\sigma})$ of 
$\fl_{H, \sigma} = {\rm Ad}_{\exp H} \fl_{\sigma}$ is
given by
\[
N_K(\fl_{H,\sigma}) \, = \, N^{'}(\fl_{H, \sigma}) K_{H, \sigma}
\, = \, K_{H, \sigma} N^{'}(\fl_{H, \sigma}).
\]
\end{prop}

\noindent
{\bf Proof.} 
It is clear from Lemma \ref{lem_lk} that 
$N^{'}(\fl_{H, \sigma})$ normalizes $\flhs$, so it 
normalizes $\flhs \cap \fk$ and the corresponding connected group
$K_{H,\sigma}$. 
This implies the second equality, and the inclusion
$K_{H, \sigma} N^{'}(\fl_{H, \sigma}) \subset 
N_K(\flhs)$.

Conversely, suppose that $k \in K$ normalizes $\flhs$. Then it
normalizes the group $K_{H,\sigma}$, so ${\rm Ad}_kT$
is a maximal torus of $K_{H,\sigma}$, where $T$ is the maximal 
torus of $K$ with Lie algebra $\ft$.
Thus there exists
$k_1 \in K_{H,\sigma}$ such that ${\rm Ad}_{k_{1}^{-1}} 
{\rm Ad}_k T = T$, i.e., 
$k_{1}^{-1}k \in N_K(T) = N_K(\ft)$.
Write $n = k_{1}^{-1}k$, so that $k = k_1 n$. It
remains to show that $n \in N^{'} (\flhs)$.

Denote by $w_n$ the Weyl group element $nT \in W$. Since $n$ normalizes
$\flhs$, it normalizes its nilradical $\fn_\sigma$. Thus $w_n \in W_\sigma$.
Now for each $\alpha \in [S_\sigma]$, the support of $\sigma$, 
consider the space
\[
V_\alpha \, = \, \flhs \cap (\fg_{\alpha} \oplus \fg_{-\alpha}).
\]
By the description of the basis of $\fl_\sigma$, we know that the Killing form
of $\fg$ restricted to $V_\alpha$ is either negative definite or
positive definite depending on whether $\sigma(\a) = 1$ or
$\sigma(\a) = -1$. Now since $n$ normalizes $\flhs$, it
permutes the spaces $V_\alpha$, for $\a \in [S_\sigma]$. But $n$ preserves
the Killing form, so $\sigma(\alpha) = 1$ implies $\sigma(w_n \alpha) = 1$. 
In other words, $w_n \Sigma_\sigma = \Sigma_\sigma$. It also
follows that $n$ normalizes $\fl_\sigma$. Therefore we have
\[
{\rm Ad}_{\exp(w_nH)} \fl _\sigma \, = \, 
{\rm Ad}_{\exp H} \fl_\sigma.
\]
An easy calculation shows that this implies $\a(H) = \a(w_nH)$
for all $\a \in [S_\sigma]$. Since $H \in \fa \cap \fm_{\sigma, 1}$
and $w_n \in W_\sigma$, it follows that $H = w_n H$. Therefore $w_n 
\in W_{H, \sigma}$, or, equivalently, $n \in N^{'}(\flhs)$.
\qed

\begin{cor}
\label{cor_quo}
Let the notation be as in Notation \ref{nota_sigma}. Then
\[
N_K(\flhs) / K_{H, \sigma} \, \cong \, N^{'}(\flhs) 
/ N^{'}(\flhs) \cap K_{H, \sigma}.
\]
\end{cor}

\begin{rem}
\label{rem_sigma}
{\em
For an extended signature $\sigma$, the group
\[
W_{0, \sigma} \, = \, \{w \in W_\sigma: w \Sigma_\sigma = \Sigma_\sigma \}
\]
contains the subgroup $R_\sigma$ generated by reflections $\{s_{\alpha}
\}$ for $\a \in \Sigma_\sigma$ as a normal subgroup.   Indeed, this
follows from the formula for $s_\alpha$ and Formula (\ref{eq_userho})
for $\sigma$. 
Set $Z_\sigma = W_{0, \sigma} / R_\sigma$. Regard $\sigma$
as a signature for the root system $[S_\sigma]$. Then $\sigma$
defines a signature for each irreducible
subsystem of $[S_\sigma]$, and we can calculate $Z_\sigma$
separately for each irreducible subsystem. The group $Z_\sigma$ is computed 
for each 
simple Lie algebra in \cite{os:compact}, Table 3,
p. 80, and explicit elements are given.
 For example, when $\fg=\fsl (n,\C)$, then if $\flstriv\not\cong
\fs\fu (n/2,n/2)$, then $Z_\sigma$ is trivial, and if $\flstriv\cong
\fs\fu (n/2,n/2)$ then $Z_\sigma$ is a group with two elements.
$Z_\sigma$ has no more than two elements except in the case when
$\fg=\fs\fo (4n,\C)$ and $\flstriv\cong \fs\fo(2n,2n)$, when
$Z_\sigma$ is the Klein 4-group. In particular, the group $W_{0, \sigma}$
can be calculated explicitly in each case. It follows that we can compute
the group $W_{H,\sigma}$ explicitly.
}
\end{rem}

\subsection{$(K, \pi_K)$-homogeneous Poisson structures on $K/T$}
\label{sec_on-KT}

In this section, we  determine all
$(K, \pi_K)$-homogeneous Poisson structures on the full
flag variety $K/T$, where $T$ is the maximal torus of $K$ with Lie 
algebra $\ft$. 

By Proposition \ref{prop_orbits-in-lagr-zero}, we only need to identify
those $\fl \in \Lagr_0$ such that $\fl \cap \fk = \ft$.
We can assume $\fl=\flhs = {\rm Ad}_{\exp H}
\fl_\sigma$, where $\sigma$ is an extended signature and 
$H \in \fa \cap \fm_{\sigma, 1}$, because the Poisson structure
on $K/T$ determined by any $\fl = {\rm Ad_k} \flhs$ for some $k \in K$
(such that $\fl \cap \fk = \ft$) will be $K$-equivariantly isomorphic to 
the one determined by $\flhs$.

\begin{prop}
\label{prop_KThomog}
Let  $\sigma$ be an extended signature and let 
$H \in \fa \cap \fm_{\sigma, 1}$.
Let $\flhs = {\rm Ad}_{\exp H} \fl_{\sigma}$. Then
$\flhs \cap \fk = \ft$ if and only if 
$\alpha(H)\not= 0$ for all $\alpha \in  \Sigma_\sigma$.
\end{prop}

\noindent
{\bf Proof.} This is a direct consequence of Lemma \ref{lem_lk}.
\qed

For every $\flhs$ such that $\flhs \cap \fk = \ft$,
denote by $\pi_{H, \sigma}$
the associated $(K, \pi_K)$-homogeneous Poisson structure on $K/T$.

\begin{cor}
\label{cor_all-on-KT}
The collection $\{\pi_{H, \sigma}\}$, as $\sigma$ runs over
all extended signatures and as $H$ takes all elements in $\fa \cap 
\fm_{\sigma, 1}$ such that $\alpha(H) \neq 0$ when $\sigma(\alpha) = 1$,
gives all $(K, \pi_K)$-homogeneous Poisson structure on $K/T$.
\end{cor}
 
An explicit formula for $\pi_{H, \sigma}$ is given in 
\cite{lu:cdyb} as 
\[
\pi_{H, \sigma} \, = \, p_* \pi_K \, + \,
{1 \over 2} \left(\sum_{\alpha \in [S_\sigma]\cap \Sigma_{+}}
{1 \over 1-\sigma(\alpha)e^{2\alpha(H)}} 
\Xa \wedge \Ya \right)^0,
\]
where $p: K \rightarrow K/T$ is the natural projection,
and the second term on the right hand side is the
$K$-invariant bi-vector field on $K/T$
whose value at $e = eT$ is the expression given
in the parenthesis. 
The fact that these are
 all the $(K, \pi_K)$-homogeneous
Poisson structures on $K/T$ up to  $K$-equivariant isomorphisms
is also proved in \cite{lu:cdyb} by a different method. Namely,
we show in \cite{lu:cdyb} that every 
such 
Poisson structure comes from a solution  to the
{\it Classical Dynamical Yang-Baxter Equation} \cite{e-v:cdyb}.
In \cite{lu:cdyb}, we also study some geometrical properties 
of these  Poisson structures such as their symplectic leaves,
modular vector fields, and moment maps for the $T$-action.

Recall from Proposition \ref{prop_normalizer} and 
Notation \ref{nota_sigma} that when
$\flhs \cap \fk = \ft$, the normalizer
subgroup $N_K(\flhs)$ of $\flhs$ in $K$ lies in the
normalizer subgroup of $\ft$ in $K$, and  we have
\[
N_K(\flhs) / T \, = \, 
W_{H, \sigma} \, = \, 
\{w \in W_\sigma: w \Sigma_\sigma = \Sigma_\sigma,
wH = H\}.
\]
When $W_{H, \sigma}$ is trivial, the Poisson manifold
$(K/T, \pi_{H, \sigma})$ embeds into $(\Lagr_0, \Pi)$
as a Poisson submanifold. When
$W_{H, \sigma}$ is not trivial, it
 follows from
Proposition \ref{prop_cover} that action of $W_{H, \sigma}$ on
$K/T$ from the right  defined by
\[
(K/T) \times W_{H, \sigma} \lrw K/T: \, (kT, w) \Map kwT
\]
is by Poisson isomorphisms. Thus, 
the group $W_{H, \sigma}$ gives symmetries of the Poisson
structure.
As we mentioned in Remark \ref{rem_sigma}, this group can 
be calculated case by case.

\begin{rem}
\label{rem_normzsigma}
{\em
If $H\in \fa$ is regular
 in the sense that it is not fixed by
any Weyl group element, then $W_{H, \sigma}$ is trivial
for any $\sigma$. 
 On the other hand, Borel and de Siebenthal showed that
every nontrivial 
signature $\sigma$ corresponding to the trivial diagram automorphism
can be put in a form such that $\sigma(\alpha_k)=-1$ for exactly one
simple root $\alpha_k$ \cite{bodeS} or \cite{os:compact}, Appendix.
 In particular, the group $W_{0,\sigma}$ contains
the Weyl group of a maximal Levi subgroup, so for
$W_{H, \sigma}$ to be trivial,
$H$ cannot be fixed by any element in a
maximal Levi subgroup, so in particular, $H$ can lie in at most
one wall.
}
\end{rem}

\begin{exam}
\label{exam_sl2-final}
{\em 
We can compute the Poisson structure $\Pi$ on $\Lagr_0$
explicitly for the case of $\fg = \fsl(2, \C)$. In this case,
it follows from \cite{dp:compactification} that $\Lagr_0$
can be $G = PSL(2, \C)$-equivariantly identified with 
$\R P^3$, regarded as the projectivization of the space
${\cal H}$ of $2 \times 2$ Hermitian matrices, where the action of
$G$ on ${\cal H}$ is by
\[
g \circ X \, = \, gX\bar{g}^t, \hspace{.2in} g \in G, X \in 
{\cal H}.
\]
The $R$-matrix $R \in \fg \wedge \fg$  (see Section \ref{sec_Lagr-general})
is explicitly given by
\[
R \, = \, -{1 \over 2} \left(
ih \wedge h - \Xa \wedge i\ea + \Ya \wedge \ea \right),
\]
where
\[
h = {1 \over 2\sqrt{2}} \left(\begin{array}{ll}1 & 0 \\ 0 & -1 \end{array}
\right), \hspace{.1in} 
\Xa = {1 \over 2} \left(\begin{array}{ll}0  & 1 \\ -1 & 0 \end{array}
\right), \hspace{.1in}
\Ya = {1 \over 2} \left(\begin{array}{ll}0 & i \\ i & 0 \end{array}
\right), \hspace{.1in}
\]
and $\ea = {1 \over 2} (\Xa -i\Ya)$. Denote by 
$v: \fg \rightarrow \chi^1({\cal H})$ the Lie algebra anti-homomorphism
defined by the above action of $G$ on ${\cal H}$, where
$\chi^1({\cal H})$ is the space of vector fields on ${\cal H}$. Then
$\Pi = {1 \over 2} v(R)$ is a Poisson structure on ${\cal H}$.
Write an element of ${\cal H}$ as
\[
X \, = \, \left(\begin{array}{ll} x & u+iv \\u-iv & y \end{array}
\right)
\]
with $x, y, u, v \in \R$. Then the Poisson brackets for $\Pi$ are
given by
\beqa
& & \{x, \, y\} \, = \, 0, \hspace{.2in}
\{x, \, u\} \, = \, -{1 \over 4}y v, \hspace{.2in}
\{x, \, v\} \, = \, {1 \over 4}y u\\
& & \{y, \, u\} \, = \, {1 \over 4}yv, \hspace{.2in}
\{y, \, v\} \, = \, -{1 \over 4}yu, \hspace{.2in}
\{u, \, v\} \, = \, {1 \over 8} y (y-x).
\eeqa
Note that
\[
c_1 \, = \, x \, + \, y \hspace{.2in} {\rm and} \hspace{.2in}
c_2 \, = \, xy \, - \, u^2 \, - \, v^2
\]
are two Casimir functions. Hence all $SU(2)$-orbits
are Poisson submanifolds.
Since this Poisson structure is quadratic, it gives rise to 
one on $\R P^3$, which is the Poisson structure $\Pi$ on $\Lagr_0$.
It can be checked that  by looking at the $SU(2)$-orbits through the 
points in $\R P^3$ corresponding to
\[
\left(\begin{array}{ll}b & 0 \\0 & 1 \end{array}
\right), \, b \in {\Bbb R}, b \neq 1
\]
we get all the $(K, \pi_K)$-homogeneous Poisson
structures $\pi_{H, \sigma}$ on $SU(2)/S^1$, 
up to $K$-equivariant isomorphisms,
as discussed in Section \ref{sec_on-KT}.  By identifying
$SU(2)/S^1$ with $S^2 =\{(x, y, z) \in \R^3: x^2 + y^2 + z^2 = 1\}$,
these Poisson structures are given by
\[
\{x, y\} = {1 \over 4} (x+2a-1) z, \hspace{.2in}
\{y, z\} = {1 \over 4} (x+2a-1) x, \hspace{.2in}
\{z, x\} = {1 \over 4} (x+2a-1)y,
\]
for $a \in \R$. Note that the antipodal map is a symmetry for the case
when $a = {1 \over 2}$. This corresponds to the fact that the
stabilizer subgroup in $SU(2)$ of the point in $\R P^3$ corresponding to
 $\left(\begin{array}{ll}1 & 0 \\ 0 & -1 \end{array}
\right)$ has two connected components.
}
\end{exam}

\end{document}